\documentclass[letterpaper,journal]{IEEEtran}
\usepackage{lineno}

\usepackage{cite}

\ifCLASSINFOpdf
\else
\usepackage[dvips]{graphicx}
\fi
\usepackage{amsmath,amsfonts}
\usepackage{array}

\ifCLASSOPTIONcompsoc
	\usepackage[caption=false,font=normalsize,labelfont=sf,textfont=sf]{subfig}
\else
	\usepackage[caption=false,font=footnotesize]{subfig}
\usepackage{dsfont,enumitem}

\newcommand{\parenths}[1]{\left( #1 \right)}
\newcommand{\brackets}[1]{\left[ #1 \right]}
\newcommand{\braces}[1]{\left\{ #1 \right\}}
\newcommand{\abs}[1]{\left| #1 \right|}
\newcommand{\norm}[1]{\left\Vert #1 \right\Vert}

\newcommand{\cond}[1]{\Bigl\lvert_{#1}}

\newtheorem{condition}{Condition}

\newtheorem{theorem}{Theorem}
\newtheorem{proposition}{Proposition}
\newtheorem{lemma}{Lemma}

\newenvironment{psmallmatrix}
  {\left(\begin{smallmatrix}}
  {\end{smallmatrix}\right)}

\begin{document}

\title{Online Maximum Likelihood Estimation of the Parameters of Partially Observed Diffusion Processes}
\author{Simone~Carlo~Surace
        and~Jean-Pascal~Pfister~\IEEEmembership{}
\thanks{Simone Carlo Surace and Jean-Pascal Pfister are affiliated to the Institute of Neuroinformatics and Neuroscience Center Zurich, University of Zurich and ETH Zurich, 8057 Zurich, Switzerland. Email: $\mathtt{surace@ini.uzh.ch}$ and $\mathtt{jpfister@ini.uzh.ch}$.}
\thanks{This work was supported by the Swiss National Science Foundation, grants PP00P3\_179060 and PP00P3\_150637.}
\thanks{This paper is an expanded version of work that was first presented at the 2016 Time Series Workshop at NIPS (Neural Information Processing Systems) conference in Barcelona, Spain.}}

\markboth{IEEE TRANSACTIONS ON AUTOMATIC CONTROL,~Vol.~XX, No.~X, XXX~XXX}%
{Shell \MakeLowercase{\textit{et al.}}: Bare Demo of IEEEtran.cls for IEEE Journals}

\maketitle

\begin{abstract}
We revisit the problem of estimating the parameters of a partially observed diffusion process, consisting of a hidden state process and an observed process, with a continuous time parameter.
The estimation is to be done online, i.e. the parameter estimate should be updated recursively based on the observation filtration.  
We provide a theoretical analysis of the stochastic gradient ascent algorithm on the incomplete-data log-likelihood.
The convergence of the algorithm is proved under suitable conditions regarding the ergodicity of the process consisting of state, filter, and tangent filter. 
Additionally, our parameter estimation is shown numerically to have the potential of improving suboptimal filters, and can be applied even when the system is not identifiable due to parameter redundancies.
Online parameter estimation is a challenging problem that is ubiquitous in fields such as robotics, neuroscience, or finance in order to design adaptive filters and optimal controllers for unknown or changing systems.
Despite this, theoretical analysis of convergence is currently lacking for most of these algorithms.
This article sheds new light on the theory of convergence in continuous time.
\end{abstract}

\section{Introduction}
\IEEEPARstart{W}{e} consider the following family of partially observed dimensional diffusion process under the probability measure $P_{\theta}$:
\begin{align}
\label{eqn:hidden} dX_t&=f(X_t,\theta)dt+g(X_t,\theta)dW_t, \\
\label{eqn:obs} dY_t&=h(X_t,\theta)dt+dV_t,
\end{align}
parametrized by $\theta\in\Theta$, where $\Theta\subset\mathds{R}^p$ is an open subset.
The process $X_t$ is called the hidden state or signal process with values in $\mathds{R}^n$, and $Y_t$ is called the observation process with values in $\mathds{R}^{n_y}$.
In addition, $W_t$, $V_t$ are independent $\mathds{R}^{n'}$- and $\mathds{R}^{n_y}$-valued standard Wiener processes (signal and observation noise).
For all $\theta\in\Theta$ we assume the initial conditions $X_0\sim p_0(\theta)$ to be independent of $W_t$ and $V_t$, we set $Y_0=0$, and we assume that $f(\cdot,\theta),g(\cdot,\theta),h(\cdot,\theta)$ are functions from $\mathds{R}^n\times\mathds{R}^p$ to $\mathds{R}^n$, $\mathds{R}^{n\times n'}$, and $\mathds{R}^{n_y}$ respectively that ensure the existence and uniqueness in probability of strong solutions to Eqs.~(\ref{eqn:hidden},\ref{eqn:obs}) for all $t\geq 0$.
Additional regularity conditions for $f,g,h$ in both arguments will be required for the convergence proof.

This setting is familiar in classical filtering theory, where the problem is to find (assuming the knowledge of $\theta$) the conditional distribution of $X_t$ conditioned on the history of observations $\mathcal{F}_{t}^{Y}=\sigma\{Y_s, 0\leq s\leq t\}$. 
In this paper, we focus on the following parameter estimation problem: 
assuming that a system with parameter $\theta_0$ generates observations $Y_t$, we want to estimate $\theta_0$ from $\mathcal{F}_{t}^{Y}$ recursively.

We will consider a well-known algorithm for parameter estimation, the so-called stochastic gradient ascent (SGA) on the incomplete-data log-likelihood function.
The stochasticity comes from the online estimation from the stream of observations, which provides a noisy estimate of the gradient of the asymptotic log-likelihood.
The main open issue we are addressing in this article is the analysis of the convergence of the parameter estimate.

This paper is structured in the following way. 
In the next section, we describe the method of obtaining recursive parameter estimates.
Next, in Section~\ref{sec:conv}, we prove the almost sure convergence of the recursive parameter estimates to stationary points of the asymptotic likelihood.
In Section~\ref{sec:ex} we provide a few numerical examples, including cases where the model is not identifiable and the filter is suboptimal.
Finally, in Section~\ref{sec:rel} we discuss the theoretical similarities and differences to related methods of recursive parameter estimation.

\section{Methods}
\label{sec:methods}
In this paper, we consider the problem of finding an estimator $\tilde\theta_t$ that is $\mathcal{F}_{t}^{Y}$-measurable and recursively computable, such as to estimate $\theta_0$ online from the continuous stream of observations.
For this task, we propose an approach based on a modification of offline maximum likelihood estimation, and therefore need to compute the likelihood of the observations (also called incomplete-data likelihood) as a function of the model parameters. 

It is a fundamental theorem of filtering theory\footnote{For a detailed exposition of the mathematical background (such as Girsanov's theorem, changes of measure, or the filtering equation \eqref{eqn:KS}), we suggest a look at the standard literature on filtering theory, e.g. \cite{Bain:2009id}.} that the \emph{innovation process} $I_t$, defined by
\begin{equation}
I_t=Y_t-\int_0^t\hat h_s(\theta) ds, \quad  \hat h_s(\theta)= \mathds{E}_{\theta}\brackets{h(X_s,\theta)\cond{}\mathcal{F}_{s}^{Y}},
\end{equation}
is a $(P_{\theta},\mathcal{F}_{t}^{Y})$-Brownian motion. 
By applying Girsanov's theorem, we can change to a measure $\tilde P$ under which $Y_t$ is a $(\tilde P,\mathcal{F}_{t}^{Y})$-Brownian motion and thus (statistically) independent of both the hidden state $X_t$ and the parameter $\theta$. 
The change of measure has a Radon-Nikodym derivative
\begin{equation}
\frac{dP_{\theta}}{d\tilde P}\cond{\mathcal{F}_{t}^{Y}}=\exp\brackets{\int_0^t \hat h_s(\theta)\cdot dY_s-\frac{1}{2}\int_0^t \Vert\hat h_s(\theta)\Vert^2 ds},
\end{equation}
where $\cdot$ denotes the Euclidean scalar product.

Since the reference measure $\tilde P$, restricted on $\mathcal{F}^Y_t$, does not depend on $\theta$, we can express the incomplete-data log-likelihood function in terms of the optimal filter as
\begin{equation}
\mathcal{L}_t(\theta)=\log\frac{dP_{\theta}}{d\tilde P}\cond{\mathcal{F}_{t}^{Y}}=\int_0^t \hat h_s(\theta)\cdot dY_s-\frac{1}{2}\int_0^t \Vert\hat h_s(\theta)\Vert^2 ds.
\label{likelihood}
\end{equation}

\subsection{Offline algorithm}
We start by describing an offline method for parameter estimation using the log-likelihood function in Eq.~\eqref{likelihood}, which serves as a basis for the online method. 

If we were interested in offline learning, our goal would be to maximize the value of $\mathcal{L}_t(\theta)$ for fixed $t$. 
There is a number of methods to solve this optimization problem. 
Among these, a simple iterative method is the gradient ascent, where an estimate $\tilde\theta_k$ at iteration $k$ is updated according to
\begin{equation}
\tilde\theta_{k+1}=\tilde\theta_k+\gamma_k \partial^{\top}_{\theta}\mathcal{L}_t(\theta)\cond{\theta=\tilde\theta_k},
\end{equation}
where $\gamma_k>0$ is called the \emph{learning rate}, and $\partial^{\top}_{\theta}$ denotes the (Euclidean) gradient operator with respect to the parameter $\theta$. 
At each iteration, the derivative of the likelihood function has to be recomputed. 
From Eq.~\eqref{likelihood}, we obtain  
\begin{equation}
\partial_{\theta}\mathcal{L}_t(\theta)=\int_0^t \parenths{dY_s-\hat h_s(\theta) ds}^{\top}\hat h_s^{\theta}(\theta),
\end{equation}
where $\cdot^{\top}$ denotes the matrix transpose and the last factor of the integrand, denoted by 
\begin{equation}
\hat h_s^{\theta}(\theta)\doteq\partial_\theta\hat h_s(\theta),
\end{equation}
takes values in the matrices of size $n_y\times p$ and is called the \emph{filter derivative} of $h$ with respect to $\theta$.\footnote{Here and in the sequel, we use the convention that the gradient operator adds a covariant dimension to the tensor field it acts on. For example, $\partial_{\theta}\mathcal{L}_t(\theta)$ takes values that are covectors (row vectors), and the gradient of $\hat h_t(\theta)$, which has values in $\mathds{R}^{n_y}$, wrt. $\theta$, is a ($n_y\times p$)--matrix ($\mathds{R}^{n_y}\otimes\mathds{R}^{p\ast}$--tensor, where $^*$ denotes a dual space)-valued process which we denote by $\hat h^{\theta}_t(\theta)$.}

In principle, computing the quantities $\hat h_t(\theta)$ requires the solution of the Kushner-Stratonovich filtering equation
\begin{equation}
d\hat\varphi_t=(\widehat{\mathcal{A}\varphi})_tdt+\parenths{(\widehat{h\varphi})_t-\hat h_t\hat\varphi_t}\cdot\parenths{dY_t-\hat h_tdt}
\label{eqn:KS}
\end{equation}
for abritrary integrable $\varphi:\mathds{R}^n\to\mathds{R}$, where $\mathcal{A}$ is the generator of the process $X_t$.
However, exact solutions are rarely available.
In the following, we assume that Eq.~\eqref{eqn:KS} admits a finite-dimensional recursive solution or a finite-dimensional recursive approximation.
This means that there is an $\mathcal{F}_{t}^{Y}$-adapted process $M_t(\theta)$ with values in $\mathds{R}^m$ and a mapping $\psi_h:\mathds{R}^m\times\Theta\rightarrow\mathds{R}^{n_y}$ such that either $\hat h_t(\theta)=\psi_h(M_t(\theta),\theta)$ (in the case of an exact solution), or such that the equation holds approximately, i.e. with some bounds (preferably uniform in time) on $$\text{Var}\|\hat h_t(\theta)-\psi_h(M_t(\theta),\theta)\|.$$
For example, in the linear-Gaussian case and if $X_0$ has a Gaussian distribution, the optimal filter can be represented in terms of a Gaussian distribution with mean $\mu_t$ and variance $P_t$, i.e. $m=2$, $M_t=(\mu_t,P_t)$, and for $h_{\theta}(x)=\theta x$, we have $\psi_h(M_t(\theta),\theta)=\theta\mu_t$. 
Apart from the linear-Gaussian case \cite{Kalman:1961wy} just mentioned, finite-dimensional (exact) recursive solutions only exist for a small class of systems, namely the Bene\v{s} class and its extensions \cite{BeneS:1981bn} --\nocite{BeneS:1985tf,Daum:1986hh,Daum:1987uh,Charalambous:1998ir}\cite{Charalambous:1998fl}.
Meanwhile, finite-dimensional recursive \emph{approximations} are available for a large class of systems, but the appropriate choice of approximation is a complex topic in its own right and will not be explored here.
We merely mention a few standard approximation schemes: extended and unscented Kalman filters \cite{Gustafsson:2012ji,Wan:2000gc}, projection or assumed-density filters \cite{Brigo:1998wr,Brigo:1999wn}, particle filters \cite{Doucet:2011us}, and particle filters without weights \cite{Evensen:2003io} --\nocite{Crisan:2010fm}\cite{Yang:2013hf}.

Given a finite-dimensional representation of the filter, a corresponding representation of the filter derivative may be formally defined by differentiation with respect to $\theta$:
\begin{equation}
\hat h^{\theta}_t(\theta)\simeq\partial_{\theta}\psi_h(M_t(\theta),\theta)+\partial_M\psi_h(M_t(\theta),\theta)M^{\theta}_t(\theta),
\end{equation}
where $\partial_M$ denotes the gradient wrt. the first argument of $\psi_h$ and $M^{\theta}_t(\theta)$ denotes the $(m\times p)$--matrix valued derivative of the process $M_t(\theta)$.
For the system in Eqs.~(\ref{eqn:hidden},\ref{eqn:obs}) and for a large class of exact and approximate filters, $M_t(\theta)$ solves a stochastic differential equation (SDE) of the form
\begin{multline}
dM_t(\theta) = \mathcal{R}(\theta,M_t(\theta))dt+\mathcal{S}(\theta,M_t(\theta))dY_t\\
+\mathcal{T}(\theta,M_t(\theta))dB_t,
\end{multline}
where $\mathcal{R},\mathcal{S}$, and $\mathcal{T}$ go to $\mathds{R}^m$, $\mathds{R}^{m\times n_y}$, and $\mathds{R}^{m\times m'}$ respectively, and $B_t$ is an $m'$-dimensional Brownian motion that is independent of $\mathcal{F}^{X,Y}_t$ (e.g. independent noise in particle filters).
By differentiating wrt. $\theta$, we find the corresponding SDE for $M^{\theta}_t(\theta)$
\begin{multline}
dM^{\theta}_t(\theta) = \mathcal{R}'(M_t(\theta),M^{\theta}_t(\theta),\theta)dt\\
 +\mathcal{S}'(M_t(\theta),M^{\theta}_t(\theta),\theta)dY_t \\
 +\mathcal{T}'(M_t(\theta),M^{\theta}_t(\theta),\theta)dB_t,
\end{multline}
where the tensor fields $\mathcal{R',S',T'}$ are given by
\begin{multline}
\mathcal{R}'(M_t(\theta),M^{\theta}_t(\theta),\theta) = \partial_{\theta}\mathcal{R}(M_t(\theta),\theta)\\
+ \partial_M\mathcal{R}(M_t(\theta),\theta)M^{\theta}_t(\theta),
\end{multline}
and analogously for $\mathcal{S}$ and $\mathcal{T}$.
In Section~\ref{sec:ex}, we will present examples of both exact and approximate filters for which these calculations will be made explicit.

These equations can be conveniently summarized in a single SDE
\begin{equation}
d\mathcal{X}_t(\theta)=\Phi(\mathcal{X}_t(\theta),\theta)dt+\Sigma(\mathcal{X}_t(\theta),\theta) d\mathcal{B}_t.
\end{equation}
Here, $\mathcal{X}_t(\theta)$ is a $D$-dimensional process defined by concatenating the state $X_t$, the filter representation and all the filter derivatives as follows:
\begin{multline}
\mathcal{X}_t(\theta)=\mathcal{C}(X_t,M_t(\theta),M^{\theta}_t(\theta)) \\
\doteq\Big(X_{t,1},...,X_{t,n},M_{t,1}(\theta),...,M_{t,m}(\theta), \\
M_{t,1}^{\theta_1}(\theta),...,M_{t,m}^{\theta_1}(\theta),...,M_{t,1}^{\theta_p}(\theta),...,M_{t,m}^{\theta_p}(\theta)\Big)^{\top},
\end{multline}
where $D=n+m+mp$, $\mathcal{C}:\mathds{R}^n\times\mathds{R}^m\times\mathds{R}^{m\times p}\rightarrow\mathds{R}^D$ is the concatenation, and $\mathcal{B}_t$ is the Wiener process defined by $\mathcal{B}_t=\parenths{W_t,V_t,B_t}$.

\subsection{Online algorithm}
Instead of integrating the gradient of the log-likelihood function up to time $t$, an SGA uses the integrand of the gradient of the log-likelihood (evaluated with the current parameter estimate) to update the parameter estimate online as new data is reaching the observer. 
The SDE for the SGA takes the form
\begin{equation}
d\tilde\theta_t=\begin{cases}
\gamma_t F(\mathcal{\tilde X}_t,\tilde\theta_t)dt+\gamma_tH(\mathcal{\tilde X}_t,\tilde\theta_t)^{\top}dV_t, &\tilde\theta_t\in\Theta \\
0, &\tilde\theta_t\notin\Theta
\end{cases},
\label{eqn:parameterSDE}
\end{equation}
where $\mathcal{\tilde X}_t\doteq\mathcal{C}(X_t,\tilde M_t,\tilde M^{\theta}_t)$ is a diffusion process with SDE
\begin{equation}
d\mathcal{\tilde X}_t=\Phi(\mathcal{\tilde X}_t,\tilde \theta_t)dt+\Sigma(\mathcal{\tilde X}_t,\tilde \theta_t) d\mathcal{B}_t,
\end{equation}
consisting of the state as well as the filter and filter derivatives integrated with the online parameter estimate.
The functions $l,F,H$, which go from $\mathds{R}^D\times\Theta$ to $\mathds{R}$, $\mathds{R}^p$, and $\mathds{R}^{n_y\times p}$ respectively, are defined as
\begin{align}
l(x,\theta)&=\psi_h(M,\theta)\cdot\brackets{h(X,\theta_0)-\tfrac{1}{2}\psi_h(M,\theta)}, \\
F(x,\theta)&\doteq H(x,\theta)^{\top}\brackets{h(X,\theta_0)-\psi_h(M,\theta)}, \\
H(x,\theta)&\doteq \partial_{\theta}\psi_h(M,\theta)+\partial_M\psi_h(M,\theta)M',
\end{align}
where $(X,M,M')=\mathcal{C}^{-1}(x)$ are the components of $x$. 
The function $l$ will be used later on (Eq.~\eqref{eqn:likely} and ff).

\section{Convergence analysis}
\label{sec:conv}

As in any stochastic gradient method, convergence relies on being able to control the errors of estimating the gradient.
This is usually done by assuming ergodicity of the system and applying regularity results on a Poisson equation, as in the treatments of related problems by \cite{LeGland:1997kp} (discrete-time), as well as \cite{Sirignano:2017if} (continuous-time but fully observed).
In our case, the ergodic system consists of the hidden state, the filter, and the filter derivative.
We therefore need to find the assumptions that guarantee that this system is ergodic, with appropriate regularity results.
We attack this problem in part \ref{sec:convB} below by giving conditions directly in terms of the finite-dimensional approximation.
However, this means that these conditions have to be checked on a case-by-case basis in order to obtain convergence results.
Once the question of ergodicity is settled, the remainder of the proof is very similar to the one in \cite{Sirignano:2017if}.

Besides this direct verification approach, the only hope of otherwise obtaining ergodicity seems to be via the \emph{optimal} filter.
This is due to the fact that the finite-dimensional system is usually highly degenerate, such that the standard theory which was used e.g. in \cite{Sirignano:2017if}, does not apply\footnote{The theory for elliptic diffusions \cite{Veretennikov:1997fu,Pardoux:2001we,Pardoux:2003wk} is clearly not applicable, and it is not clear how to apply hypoellipticity or H\"ormander's condition \cite{Hormander:1967ik,Bony:1969tf,Ichihara:1974ig} in general. For example, in the linear-Gaussian case the process $\mathcal{X}_t$ does not satisfy the parabolic H\"ormander condition.}.
Ergodicity of the \emph{optimal} filter for a stochastic dynamical system of the form of Eqs.~(\ref{eqn:hidden},\ref{eqn:obs}) follows from the ergodicity of the hidden state process and the non-degeneracy of observations\footnote{Note that \cite{Kunita:1971eo} contained a gap that has been fully closed by \cite{vanHandel:2009ev}. } (see \cite{Kunita:1971eo,Kunita:1991uy,Budhiraja:2003hz,vanHandel:2009ev}).
The problem is then to extend these results to the derivative of the optimal filter with respect to the parameters, and to transfer them to \emph{approximate} finite-dimensional representations of the filter, given some bounds on the accuracy of the approximation.
The question of transferring ergodicity of the \emph{exact} filter and filter derivative to the approximate ones, as well as the ergodicity of the filter derivative, remains open.

\subsection{Direct conditions for the ergodicity of the approximate filter}
\label{sec:convB}
Here, we give sufficient conditions directly in terms of the approximate filtering equation.
Before stating the conditions, we introduce the following notation: 
We say that a function $G:\mathds{R}^d\times\Theta\rightarrow\mathds{R}$ has the polynomial growth property (PGP) if there are $q,K>0$ such that for all $\theta\in\Theta$,
\begin{equation}
|G(x,\theta)|\leq K(1+||x||^q).
\end{equation}
Let $\mathds{G}^d$ be the function space defined by all functions $G:\mathds{R}^d\times\Theta\rightarrow\mathds{R}$ such that
\begin{enumerate}[label=(\alph*)]
\item $G(\cdot,\theta)\in C(\mathds{R}^d)$, 
\item $G(x,\cdot)\in C^2(\Theta)$, 
\item $\partial_{\theta}G(x,\cdot)$ and $\partial^2_{\theta}G(x,\cdot)$ are H\"older continuous with exponent $\alpha>0$.
\end{enumerate}
Let $\mathds{G}^d_c$ be the subset consisting of all $G\in\mathds{G}^d$ that are centered, i.e. $\int_{\mathds{R}^d}G(x,\theta)\mu_{\theta}(dx)=0$.
Let $\mathds{\bar G}^d$ be the subset consisting of all $G\in\mathds{G}^d$ such that $G$ and all its first and second derivatives wrt. $\theta$ satisfy the PGP.

Now, we may state the conditions on the processes:

\begin{condition}
\label{condfinitedim}
\begin{enumerate}[label=(\roman*)]
\item The process $\mathcal{X}_t(\theta)$ is ergodic under $P_{\theta_0}$, with a unique invariant probability measure $\mu_{\theta}$ on $(\mathds{R}^D,\mathds{B}_D)$, where $\mathds{B}_D$ is the Borel $\sigma$-algebra on $\mathds{R}^D$.
\item For any $q>0$ and $\theta\in\Theta$ there is a constant $K_q>0$ such that
\begin{equation}
\int_{\mathds{R}^D}\parenths{1+||x||^q}\mu_{\theta}(dx)\leq K_q.
\label{eqn:lemma3a}
\end{equation}
\item Define the finite signed measures $\nu_{\theta,i}=\partial_{\theta_i}\mu_{\theta}$, $i=1,...,p$ and let $|\nu_{\theta,i}(dx)|$ be their total variation.
For any $q>0$ and $\theta\in\Theta$ there is a constant $K_q'>0$ such that
\begin{equation}
\int_{\mathds{R}^D}\parenths{1+||x||^q}|\nu_{\theta,i}(dx)|\leq K_q'.
\label{eqn:lemma3ab}
\end{equation}
\item Let $\mathcal{A_X}$ be the infinitesimal generator of $\mathcal{X}_t(\theta)$ and let $G\in\mathds{G}^D_c$. 
Then the Poisson equation $\mathcal{A_X}v(x,\theta)=G(x,\theta)$ has a unique solution $v(x,\theta)$ that lies in $\mathds{G}^D$, with $v(\cdot,\theta)\in C^2(\mathds{R}^D)$. 
Moreover, if $G\in\mathds{\bar G}^D$, then $v\in\mathds{\bar G}^D$ and also $\partial_x\partial_{\theta}v$ has the PGP.
\item For all $q>0$, $\mathds{E}[||\mathcal{\tilde X}_t||^q]<\infty$ and there is a $K>0$ such that for $t$ large enough,
\begin{align}
\forall \theta\in\Theta \quad &\mathds{E}_{\theta_0}\brackets{\sup_{s\leq t}||\mathcal{X}_s(\theta)||^q}\leq K \sqrt{t}, \\
&\mathds{E}_{\theta_0}\brackets{\sup_{s\leq t}||\mathcal{\tilde X}_s||^q}\leq K \sqrt{t}.
\end{align}
\end{enumerate}
\IEEEQEDclosed
\end{condition}

\begin{condition}
\label{condregularity}
The function $F$ is in $\mathds{\bar G}^d$ (component-wise).
The function $\psi_h$ is in $\mathds{G}^m$ and has the PGP (component-wise).
In addition, $l(x,\theta)$, $H(x,\theta)$, and $\Sigma$ have the PGP (component-wise).
\IEEEQEDclosed
\end{condition}

Lastly, the following condition on the learning rate is imposed:

\begin{condition}
\label{condlearningrate}
$
\int_0^{\infty}\gamma_t dt=\infty, \; \int_0^{\infty}\gamma_t^2 dt=0
$,
and there is an $r>0$ such that $\lim_{t\rightarrow\infty}\gamma_t^2t^{1/2+2r}=0$.
\IEEEQEDclosed
\end{condition}

\subsection{Results}
Let the approximate (in the sense of using the approximate filter representation) incomplete-data log-likelihood be given by
\begin{equation}
\mathcal{L}_t(\theta)=\int_0^t l(\mathcal{X}_s(\theta),\theta)ds+\int_0^t \psi_h(M_s(\theta),\theta)\cdot dV_s.
\label{eqn:likely}
\end{equation}
Under the above conditions 1--3, we have the following

\begin{proposition}[Regularity of the asymptotic likelihood]
\label{prop1}
\begin{enumerate}[label=(\roman*)]
\item The process $\frac{1}{t}\mathcal{L}_t(\theta)$ converges almost surely (a.s.) to $\mathcal{\tilde L}(\theta)$, which is given by
\begin{align}
\mathcal{\tilde L}(\theta) &= \int_{\mathds{R}^D}l(x,\theta)\mu_{\theta}(dx).
\end{align}
\item The asymptotic likelihood function $\mathcal{\tilde L}(\theta)$ is in $C^2(\Theta)$, and the gradient $g$ and Hessian $\mathcal{H}$ of the asymptotic likelihood are given in terms of the invariant measure $\mu_{\theta}$ and its derivative $\nu_{\theta}$ as
\begin{align}
g(\theta) &\doteq \partial_\theta\mathcal{\tilde L}(\theta)  = \int_{\mathds{R}^D}F(x,\theta)^{\top}\mu_{\theta}(dx),\\
\begin{split}
\mathcal{H}(\theta) &\doteq \partial^{\top}_\theta\partial_\theta\mathcal{\tilde L}(\theta)  = \int_{\mathds{R}^D}\partial_{\theta} F(x,\theta)\mu_{\theta}(dx) \\
&\quad +\int_{\mathds{R}^D} F(x,\theta)\nu_{\theta}(dx),\label{eqn:Hessian}
\end{split}
\end{align}
\item The function $G(x,\theta)\doteq F(x,\theta)-g(\theta)$ is in $\mathds{\bar G}^D\cap\mathds{G}^D_c$.
\item There is a constant $C>0$ such that
\begin{equation}
\mathcal{\tilde L}(\theta)+\norm{g(\theta)}+\norm{\mathcal{H}(\theta)}\leq C.
\label{eqn:lemma3b}
\end{equation}
\end{enumerate}
\end{proposition}
\begin{IEEEproof}
See Appendix~\ref{propproof2}.
\end{IEEEproof} 

Now we can formulate our main result.
Its proof relies on several lemmas that are given in Appendix~\ref{lemmas}.
\begin{theorem}[main theorem]
\label{mainresult}
Assume conditions~\ref{condfinitedim}-\ref{condlearningrate} and let $\tilde\theta_0\in\Theta$.
Then, with probability one
\begin{equation}
\lim_{t\rightarrow\infty} \norm{g(\tilde\theta_t)}=0 \quad \text{or} \quad \tilde\theta_t\rightarrow\partial\Theta.
\end{equation}
\end{theorem}
\begin{IEEEproof}
See Appendix~\ref{proof}.
\end{IEEEproof}

\section{Examples and numerical validation}
\label{sec:ex}
Here, we consider two different example filtering problems and show explicitly how the parameter learning rules are derived. 
We also study the numerical performance of the learning method.
Since under suitable conditions on the decay of the learning rate, convergence is guaranteed by the results in the preceding section, we do not study this case.
Instead, we study whether the method also converges with constant learning rate, i.e. when violating Condition~\ref{condlearningrate}. 
A constant learning rate is a sensible choice when the system parameters are expected to change.

All numerical experiments use the Euler-Maruyama method to integrate the SDEs. 
We evaluate the performance of the learned filter by the mean squared error (MSE), normalized by the variance of the hidden process. 


\subsection{One-dimensional Kalman-Bucy filter (linear filtering problem)}
\label{ex1}
We shall first consider the simple case of the linear filtering problem, for which it is possible to obtain an exact finite-dimensional filter as well as exact expressions for the asymptotic likelihood. 
Here, we have a three-dimensional parameter vector $\theta=(a,\sigma,w)$, where $a,\sigma>0$ and $w\in\mathds{R}$, and we have $f(x,\theta)=-a x$, $g(x,\theta)=\sigma$ and $h(x,\theta)=w x$, such that the filtering problem reads
\begin{align}
dX_t=-a X_t dt+\sigma dW_t, \qquad dY_t=wX_tdt+dV_t.
\end{align}

Assuming a Gaussian initialization, i.e. $X_0\sim\mathcal{N}(0,\sigma^2/2a)$, the optimal filter has a Gaussian distribution with mean $\mu_t$ and variance $P_t$ (the Kalman-Bucy filter \cite{Kalman:1961wy}).
This is a two-dimensional representation with $M_t(\theta)=(\mu_t(\theta),P_t(\theta))^{\top}$, which can be expressed as
\begin{multline}
dM_t(\theta)=\begin{pmatrix}
-a\mu_t(\theta)-w^2\mu_t(\theta)P_t(\theta)\\
\sigma^2-2a P_t(\theta)-w^2 P_t(\theta)^2
\end{pmatrix}dt\\
+\begin{pmatrix}
wP_t(\theta)\\
0
\end{pmatrix}dY_t.
\end{multline}
We have $\psi_h(M_t(\theta),\theta)=w\mu_t(\theta)$.

Let us first calculate the asymptotic log-likelihood.
It follows from the above that $P_t(\theta)$ (and its derivatives with respect to $\theta$) will tend to a unique steady state given by
\begin{equation}
P_{\infty}(\theta)=\frac{1}{w^2}\parenths{\sqrt{a^2+w^2\sigma^2}-a}.
\end{equation}
By initializing the filter with this steady-state value, the representation can be made one-dimensional, i.e.
\begin{multline}
dM_t(\theta)=\parenths{-a\mu_t(\theta)-w^2\mu_t(\theta)P_{\infty}(\theta)}dt\\
+wP_{\infty}(\theta)dY_t.
\end{multline}
The process $\mathcal{X}_t(\theta)$ consisting of $X_t$, $\mu_t(\theta)$, and the filter derivatives $\mu^a_t(\theta),\mu^{\sigma}_t(\theta),\mu^w_t(\theta)$, therefore admits the SDE representation
\begin{equation}
d\mathcal{X}_t(\theta)=A\mathcal{X}_t(\theta)dt+B \begin{pmatrix}
dW_t\\
dV_t
\end{pmatrix},
\end{equation}
with matrices
\begin{equation}
\setlength\arraycolsep{1pt}
A=\begin{psmallmatrix}
-a_0 & 0 & 0 & 0 & 0 \\
w w_0 P & -a-w^2 P & 0 & 0 & 0 \\
w w_0 P^{a} & -1 -w^2 P^{a} & -a-w^2P & 0 & 0 \\
w w_0 P^{\sigma} & -w^2 P^{\sigma} & 0 & -a-w^2P & 0 \\
w_0 (P+wP^w) & -w^2 P^{w} & 0 & 0 & -a-w^2P
\end{psmallmatrix},
\end{equation}
and 
\begin{equation}
\setlength\arraycolsep{1pt}
B=\begin{psmallmatrix}
\sigma_0 & 0 \\
0 & w P \\
0 & w P^a \\
0 & w P^{\sigma} \\
0 & P+w P^w
\end{psmallmatrix},
\end{equation}
where $P$ is a short-hand for $P_{\infty}(\theta)$ and $P^a$ etc. are partial derivatives of $P_{\infty}(\theta)$.

The process $\mathcal{X}_t(\theta)$ is ergodic, and its unique invariant probability measure is multivariate Gaussian with zero mean and covariance matrix $K$ given by the solution to
\begin{equation}
BB^{\top}+AK+KA^{\top}=0.
\end{equation}
In terms of this, the asymptotic log-likelihood reads
 \begin{equation}
 \begin{split}
\mathcal{\tilde L(\theta)} &= w w_0 K_{12}-\frac{1}{2}w^2K_{22} \\
&=\frac{P_{\infty}(\theta) w^2 \sigma_0^2w_0^2(2a+P_{\infty}(\theta)w^2)}{4a_0(a+P_{\infty}(\theta)w^2)(a+a_0+P_{\infty}(\theta)w^2)}\\
&\quad -\frac{P_{\infty}(\theta)^2 w^4}{4(a+P_{\infty}(\theta)w^2)}.
\end{split}
\label{eqn:asympllKB}
\end{equation}
With suitable boundaries of the parameter space, all the items from Condition~\ref{condfinitedim} can be verified.

This model is non-identifiable from the observations.
The set of critical points of the asymptotic likelihood is characterized by
\begin{equation}
\partial_{\theta}\mathcal{\tilde L(\theta)} = 0 \Leftrightarrow \theta=\parenths{a_0,\sigma, \frac{w_0\sigma_0}{\sigma}}^{\top}, \quad \sigma>0,
\end{equation}
i.e. convergence can be guaranteed to one of these points only, and not to the ground truth parameters $\theta_0=(a_0,\sigma_0,w_0)^{\top}$.
The model becomes identifiable if either $\sigma_0$ or $w_0$ is known. 
Alternatively, one may fix a parametrization for which $X_t$ has unit variance (i.e. $\sigma=\sqrt{2a}$). 

Let us now derive the parameter update equations.
The filtering equations for the mismatched filter, expressed in terms of the online parameter estimates, read
\begin{align}
\label{KBmean}d\mu_t&=-\tilde a_t \mu_t dt+\tilde w_tP_t(dY_t-\tilde w_t \mu_t dt), \quad \mu_0=0, \\
\label{KBvar}dP_t&=\parenths{\tilde\sigma_t^2-2\tilde a_t P_t-\tilde w_t^2P_t^2}dt, \quad P_0=\frac{\tilde\sigma_0^2}{2\tilde a_0},
\end{align}
where the initialization of $P_0$ reflects the prior belief of the variance of $X_0$ based on the initial parameter estimates.

The online parameter update equations read
\begin{align}
d\tilde a_t&=\gamma_a\tilde a_t\tilde w_t\mu_t^a \parenths{dY_t-\tilde w_t\mu_tdt}, \\
d\tilde \sigma_t&=\gamma_{\sigma}\tilde \sigma_t\tilde w_t\mu_t^{\sigma} \parenths{dY_t-\tilde w_t\mu_tdt}, \\
d\tilde w_t&=\gamma_w\tilde w_t\parenths{\mu_t+\tilde w_t\mu_t^w}  \parenths{dY_t-\tilde w_t\mu_tdt}.
\end{align}
In order to prevent sign changes of the parameters we chose time-dependent learning rates that are proportional to the parameters ($\tilde a_t$ has to stay non-negative because the filter equations turn unstable otherwise; for $\tilde\sigma_t$ and $\tilde w_t$ it is because of identifiability, i.e. the signs of $\sigma$ and $w$ are not identifiable from $\mathcal{F}_{t}^{Y}$). 
Here, we introduced the filter derivatives $\mu_t^a,\mu_t^{\sigma}$ and $\mu_t^w$ of the mean, which, together with the filter derivatives of the variance, satisfy the coupled system of SDEs 
\begin{align}
\begin{split}
\label{KBfilterderiv1}d\mu_t^a&=-\brackets{\mu_t+\parenths{\tilde a_t+\tilde w_t^2P_t}\mu_t^a+\tilde w_t^2\mu_tP_t^a}dt\\
&\quad+\tilde w_t P_t^a dY_t,
\end{split}\\
dP_t^a&=-\brackets{2P_t+2\parenths{\tilde a_t+\tilde w_t^2P_t}P_t^a}dt,\\
\begin{split}
d\mu_t^{\sigma}&=-\brackets{\parenths{\tilde a_t+\tilde w_t^2P_t}\mu_t^{\sigma}+\tilde w_t^2\mu_tP_t^{\sigma}}dt\\
&\quad+\tilde w_t P_t^{\sigma} dY_t,
\end{split}\\
dP_t^{\sigma}&=\brackets{2\tilde\sigma_t-2\parenths{\tilde a_t+\tilde w_t^2P_t}P_t^{\sigma}}dt,\\
\begin{split}
d\mu_t^{w}&=-\brackets{2\tilde w_t\mu_tP_t+\parenths{\tilde a_t+\tilde w_t^2P_t}\mu_t^w}dt\\
&\quad-\tilde w_t^2\mu_tP_t^wdt+\brackets{P_t+\tilde w_t P_t^{w}}dY_t,
\end{split}\\
\label{KBfilterderivlast}dP_t^{w}&=-\brackets{2\tilde w_tP_t^2+2\parenths{\tilde a_t+\tilde w_t^2P_t}P_t^w}dt,\\
\mu_0^a&=\mu_0^{\sigma}=\mu_0^{w}=0, \\
P_0^a&=-\frac{\tilde\sigma_0^2}{2\tilde a_0^2}, \quad P_0^{\sigma}=\frac{\tilde\sigma_0}{\tilde a_0}, \quad P_0^w=0.
\end{align}
The right-hand sides of the filter derivative equations and the initial conditions of the filter derivatives are obtained from the corresponding equations of the filtered mean and variance and their initial conditions by differentiating with respect to each of the parameters (see Section~\ref{sec:methods} for details).

\begin{figure*}[!t]
\centering
\subfloat{\includegraphics[width=\textwidth,natwidth=1.27\textwidth,natheight=10cm]{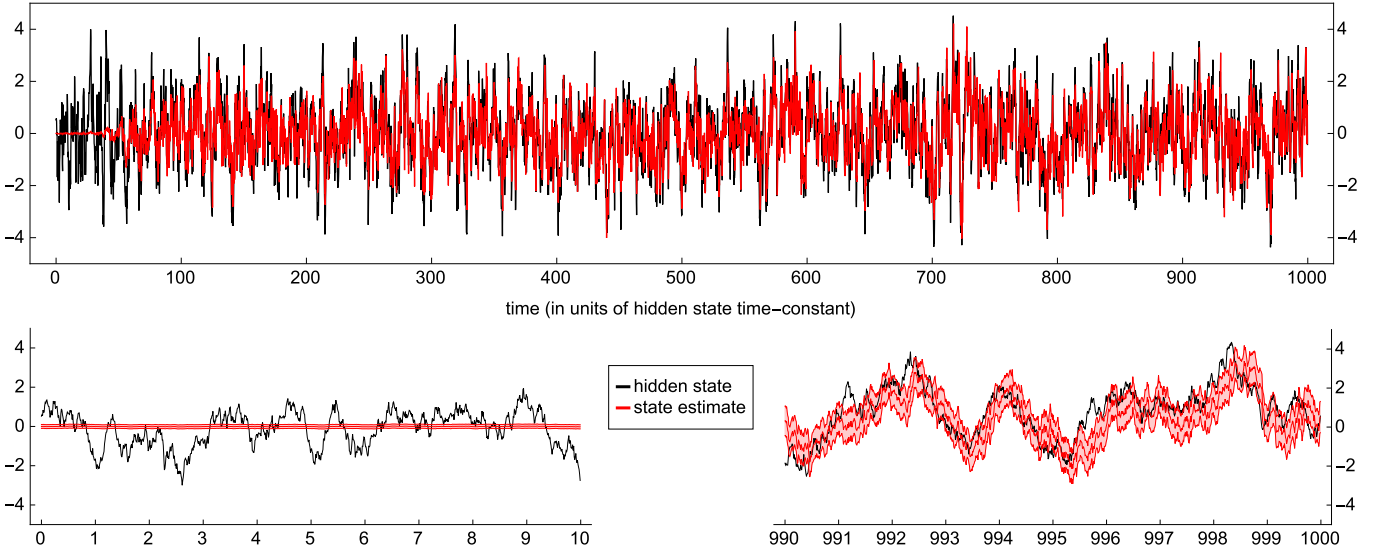}%
}
\caption{\textbf{Online learning and filtering in the linear model}. 
The hidden state $X_t$ (black) and Kalman-Bucy state estimate $\mu_t$ (red, shaded region shows $\mu_t$ $\pm$ one standard deviation $\sqrt{P_t}$, c.f. Eqs.~(\ref{KBmean},\ref{KBvar})) are shown for the linear model of Section \ref{ex1} with parameters $a_0=1$, $\sigma_0=2$, $w_0=3$. 
The time-step is $dt=10^{-3}$, initial parameter estimates are $\tilde a_0=10$, $\tilde \sigma_0=\sqrt{0.2}$, $\tilde w_0=3$ (i.e. the parameter $w_0$ is known), and the learning rates are $\gamma_a=\gamma_{\sigma}=0.03$ and $\gamma_{w}=0$. 
Top: the entire learning period of $T=1000$ shows a gradual improvement of the performance of the filter. 
Bottom left: during the first 10 seconds, the model is still strongly mismatched. 
Bottom right: during the last 10 seconds, the filter optimally tracks the hidden state.}
\label{Fig:1}
\end{figure*}

\begin{figure*}[!tb]
\centering
\subfloat{\includegraphics[width=\textwidth,natwidth=1.27\textwidth,natheight=7cm]{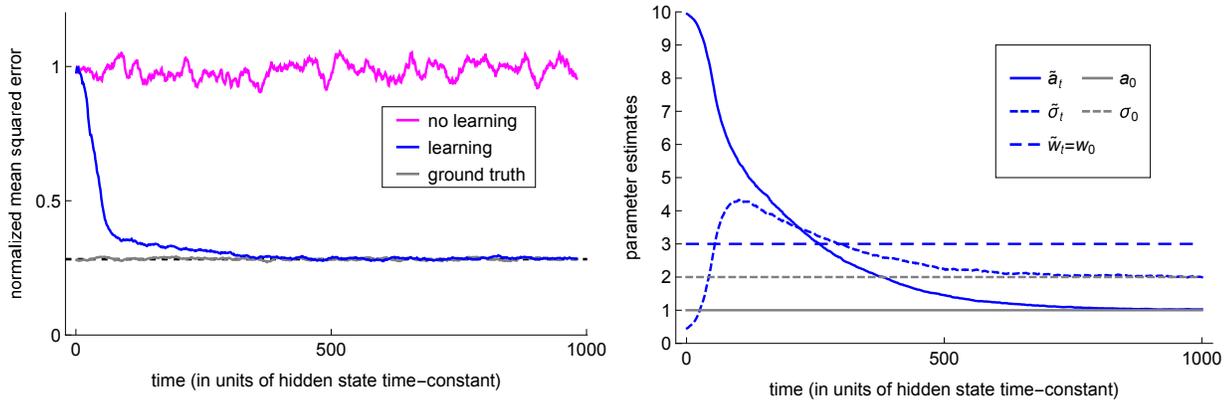}%
}
\caption{\textbf{Online learning and filtering in the linear model}. 
The time evolution of the MSE and parameter estimates are shown for the linear model of Section \ref{ex1} (see Fig.~\ref{Fig:1} caption for details). 
Left: the moving average of the normalized MSE (time window of 20 seconds) shows how the learning algorithm leads to a gradual improvement of the performance of the filter, which eventually reaches the performance of an optimal Kalman-Bucy filter with ground truth parameters. 
The black, dashed line shows the theoretical result for the performance of the Kalman-Bucy filter. 
Right: the parameter estimates for the unknown parameters converge to the ground truth parameters. 
All curves are trial-averaged ($N=100$ trials).}
\label{Fig:2}
\end{figure*}

\begin{figure*}[!tb]
\centering
\subfloat{\includegraphics[width=\textwidth,natwidth=1.3\textwidth,natheight=7cm]{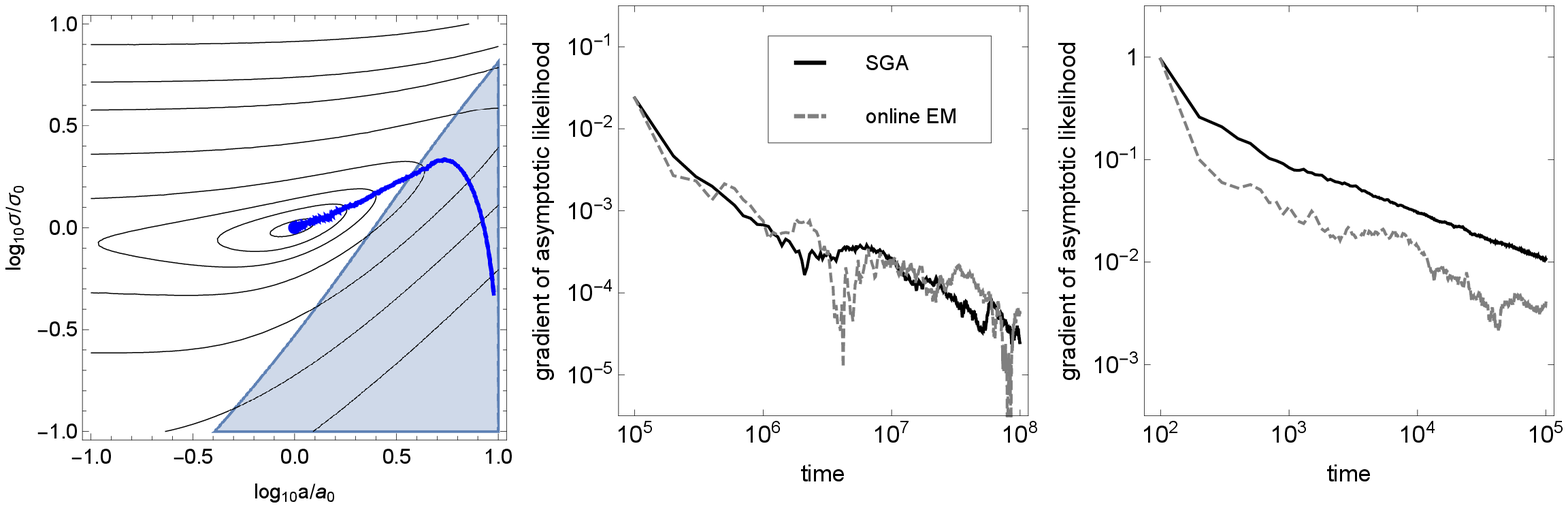}%
}
\caption{\textbf{Parameter estimation in the linear model}. 
Left: the asymptotic log-likelihood function from Eq.~\eqref{eqn:asympllKB} in the parameter subspace spanned by $a$ and $\sigma$ for $w=w_0=3$ has a single global maximum near $a=a_0$ and $\sigma=\sigma_0$. 
The shading shows the region where the function is non-concave, and the blue line is the trial-averaged learning trajectory from Fig.~\ref{Fig:2}.
Center: Single-trial convergence of the absolute value of the gradient of the asymptotic log-likelihood evaluated at the online parameter estimate given by stochastic gradient ascent (SGA) and online EM, respectively.
The rate of convergence is similar for both algorithms.
Right: An average over 10 samples reveals that the online EM is initially faster and more variable, but the order of convergence is the same. 
Time is measured in units of the intrinsic time-constant $1/a_0$.}
\label{Fig:3}
\end{figure*}

First, we investigated one of the cases where the model is identifiable, i.e. the parameter $w$ was assumed to be known and we set $\tilde w_0=w_0=3$ and $\gamma_w=0$. 
The performance of the algorithm is visualized in Fig.~\ref{Fig:1} where the learning process is shown in a single trial, and in Fig.~\ref{Fig:2}, where we show trial-averaged learning curves for the MSE and the parameter estimates. 
For both figures, the ground truth parameters were set to $a_0=1$, $\sigma_0=2$, and the initial parameter estimates were $\tilde a_0=10$ and $\tilde\sigma_0=\sqrt{0.2}$, making for a strongly mismatched model that produces an MSE close to 1 without learning, i.e. with all learning rates set to zero. 
With constant learning rates $\gamma_a=\gamma_{\sigma}=0.03$, the filter performance can be improved to almost optimal performance within a time-frame of $T=1000$, after which the parameter estimates approach the ground truth. 
The log-likelihood function is not globally concave, but it has a single global maximum (see Fig.~\ref{Fig:3} left).

The comparison to online expectation maximization (see Section \ref{sec:oEM}) is shown in Fig.~\ref{Fig:3} center and right.
Here, only the parameter $a$ is learned, while $\sigma=\sigma_0=2$ and $w=w_0=3$.
The simulations suggest that online is slightly faster in the beginning, but the order of convergence is similar for SGA and online EM.
This goes along with very similar computational complexity: the number of stochastic differential equations that have to be integrated is the same for SGA and online EM.

\subsection{Bimodal state and linear observation model with (approximate) projection filter}
\label{ex2}
\begin{figure*}[!t]
\centering
\subfloat{\includegraphics[width=\textwidth,natwidth=1.27\textwidth,natheight=10cm]{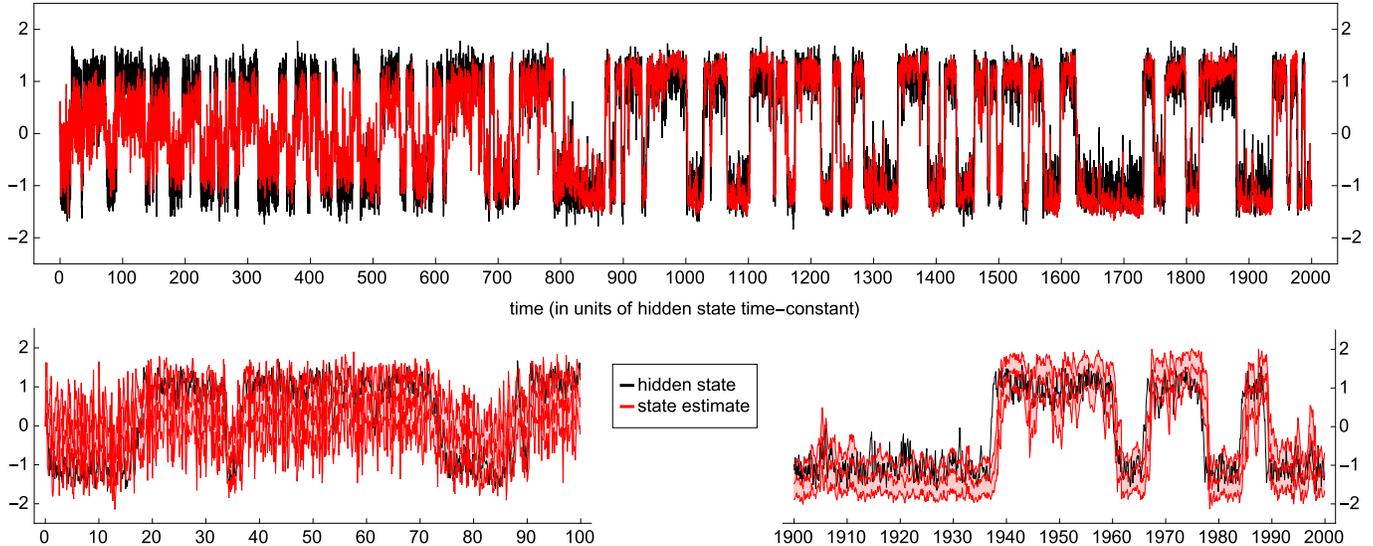}%
}
\caption{\textbf{Online learning and filtering in the nonlinear model}. 
The hidden state $X_t$ (black) and mean $\mu_t$ of the projection filter are shown for the bimodal model of Section \ref{ex2} with parameters $a_0=4$, $b_0=3$, $\sigma_0=1$ and $w_0=2$, $\tilde a_0=1$, $\tilde b_0=2$, $\tilde \sigma_0=3$, $\tilde w_0=4$, $\gamma_a=\gamma_b=\gamma_{w}=10^{-1}$ and $\gamma_{\sigma}=0.04$. 
Top: the entire learning period of $T=2000$ shows an improvement in both step size between the two attractors and the variability within both attractors. 
Bottom left: during the first 100 seconds, the filter is too sensitive to observations and has an incorrect spacing between attractors. 
Bottom right: during the last 100 seconds, the filter shows good tracking performance.}
\label{Fig:4}
\end{figure*}

Consider the following system with four positive parameters $(a,b,\sigma,w)$:
\begin{align}
dX_t&=X_t\parenths{a-b X_t^2} dt+\sigma dW_t,\\
dY_t&=wX_tdt+dV_t.
\end{align}
In this problem the hidden state $X_t$ has a bimodal stationary distribution with modes at $x=\pm\sqrt{a/b}$. 
Since the observation model is linear like in Section \ref{ex1}, the parameter learning rules are expressed in terms of the posterior mean $\mu_t=\hat X_t$ as
\begin{align}
d\tilde a_t&=\gamma_a\tilde a_t\tilde w_t\mu_t^a \parenths{dY_t-\tilde w_t\mu_tdt}, \\
d\tilde b_t&=\gamma_b\tilde b_t\tilde w_t\mu_t^b \parenths{dY_t-\tilde w_t\mu_tdt}, \\
d\tilde \sigma_t&=\gamma_{\sigma}\tilde \sigma_t\tilde w_t\mu_t^{\sigma} \parenths{dY_t-\tilde w_t\mu_tdt}, \\
d\tilde w_t&=\gamma_w\tilde w_t\parenths{\mu_t+\tilde w_t\mu_t^w}  \parenths{dY_t-\tilde w_t\mu_tdt},
\end{align}
We have made the learning rules proportional to the parameters in order to prevent sign changes, i.e. to guarantee that all parameters remain positive. 
In contrast to the linear model in Section \ref{ex1}, the filtering problem is not exactly solvable. 
We use the projection filter on the manifold of Gaussian densities introduced by \cite{Brigo:1999wn}, or equivalently, the Gaussian assumed density filter (ADF) in Stratonovich calculus. 
The mean $\mu_t$ and variance $P_t$ of the Gaussian approximation to the filter evolve as
\begin{align}
\label{bimodalmean}
\begin{split} 
d\mu_t&=\brackets{\tilde a_t\mu_t-\tilde b_t\mu_t^3-\parenths{3\tilde b_t+\tilde w_t^2}\mu_tP_t}dt\\
&\quad+\tilde w_t P_t dY_t, \quad \mu_0=0,
\end{split}\\
\label{bimodalvar}dP_t&=\brackets{\tilde\sigma_t^2+\parenths{2\tilde a_t-\tilde w_t^2 P_t^2-6\tilde b_t(\mu_t^2+P_t)}P_t}dt,
\end{align}
where the initial variance as a function of the initial parameter estimates is the variance of the stationary distribution obtained by solving the equation $\mathcal{A}^{\dag}=0$:
\begin{align}
P_0=\Gamma\parenths{\tilde a_0,\tilde b_0,\tilde\sigma_0}=\frac{\int_{-\infty}^{\infty}x^2e^{\tilde\sigma_0^{-2}\parenths{\tilde a_0 x^2-\tfrac{1}{2}\tilde b_0 x^4}}dx}{\int_{-\infty}^{\infty}e^{\tilde\sigma_0^{-2}\parenths{\tilde a_0 x^2-\tfrac{1}{2}\tilde b_0 x^4}}dx}.
\end{align}
By differentiating Eqs.~(\ref{bimodalmean},\ref{bimodalvar}) with respect to the parameters, we obtain the following equations for the filter derivatives:
\begin{align}
d\mu_t^{a}&=\brackets{\mu_t+\alpha_t\mu_t^{a}+\beta_tP_t^{a}}dt+\tilde w_tP_t^{a}dY_t,\\
dP_t^{a}&=\brackets{2P_t+A_t\mu_t^{a}+B_tP_t^{a}}dt, \\
\begin{split}
d\mu_t^{b}&=\brackets{-\mu_t\parenths{\mu_t^2+3P_t}+\alpha_t\mu_t^{b}+\beta_tP_t^{b}}dt\\
&\quad+\tilde w_tP_t^{a}dY_t,
\end{split}\\
dP_t^{b}&=\brackets{-6P_t\parenths{\mu_t^2+P_t}+A_t\mu_t^{b}+B_tP_t^{b}}dt, \\
d\mu_t^{\sigma}&=\brackets{\alpha_t\mu_t^{\sigma}+\beta_tP_t^{\sigma}}dt+\tilde w_tP_t^{\sigma}dY_t,\\
dP_t^{\sigma}&=\brackets{2\tilde\sigma_t+A_t\mu_t^{\sigma}+B_tP_t^{\sigma}}dt, \\
\begin{split}
d\mu_t^{w}&=\brackets{-2\tilde w_t\mu_tP_t+\alpha_t\mu_t^{w}+\beta_tP_t^{w}}dt\\
&\quad +\brackets{P_t+\tilde w_tP_t^{w}}dY_t,
\end{split}\\
dP_t^{w}&=\brackets{-2\tilde w_tP_t^2+A_t\mu_t^{w}+B_tP_t^{w}}dt, \\
\mu_0^a&=\mu_0^b=\mu_0^{\sigma}=\mu_0^{w}=0, \\
P_0^a&=\frac{\partial}{\partial\tilde a_0}\Gamma\parenths{\tilde a_0,\tilde b_0,\tilde\sigma_0}, \quad P_0^b=\frac{\partial}{\partial\tilde b_0}\Gamma\parenths{\tilde a_0,\tilde b_0,\tilde\sigma_0}, \\ 
P_0^{\sigma}&=\frac{\partial}{\partial\tilde \sigma_0}\Gamma\parenths{\tilde a_0,\tilde b_0,\tilde\sigma_0}, \quad P_0^w=0,
\end{align}
where we introduced the following auxiliary processes
\begin{align}
\alpha_t&=\tilde a_t-\tilde w_t^2P_t-3\tilde b_t\parenths{\mu_t^2+P_t},\\
\beta_t&=-\parenths{\tilde w_t^2+3\tilde b_t}\mu_t,\\
A_t&=-12\tilde b_t\mu_tP_t,\\
B_t&=2\tilde a_t-2\tilde w_t^2P_t-6\tilde b_t\parenths{\mu_t^2+2P_t}.
\end{align}
We numerically tested the learning algorithm for this nonlinear model by simulating a system with $a_0=4$, $b_0=3$, $\sigma_0=1$ and $w_0=2$, leading to a variance $\text{Var}(X_t)=1.17$. 
Initial parameter estimates were set to a permutation of the ground truth, i.e. $\tilde a_0=1$, $\tilde b_0=2$, $\tilde\sigma_0=3$ and $\tilde w_0=4$ and the simulations lasted $T=2000$ (due to the longer time-scale compared to the linear model) with a time-step of $dt=10^{-3}$. In Fig.~\ref{Fig:4} we show an example of the learning process.

\begin{figure*}[!t]
\centering
\subfloat{\includegraphics[width=\textwidth,natwidth=1.27\textwidth,natheight=7cm]{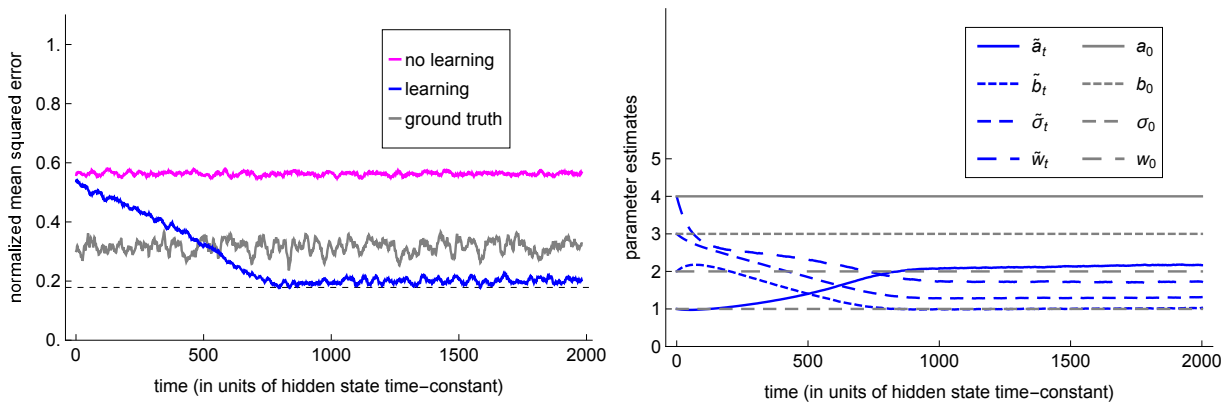}%
}
\caption{\textbf{Online learning and filtering in the nonlinear model}. 
The time evolution of the MSE and parameter estimates are shown for the bimodal model of Section \ref{ex2} (see Fig.~\ref{Fig:4} caption for details). 
Left: the moving average of the normalized MSE (time window of 20 seconds) shows how the learning algorithm allows the filter performance to improve to a level that is better than that of a filter with fixed parameters set to the ground truth. 
However, it is still slightly worse than an optimal filter; the dashed black line shows the performance of a particle filter with 1000 particles with parameters set to the ground truth. 
Right: despite the low filter error, the parameter estimates do \emph{not} converge to the ground truth. 
All curves are trial-averaged ($N=100$ trials).}
\label{Fig:5}
\end{figure*}

In this case, the sub-optimality of the Gaussian approximation inherent in the projection filter allows the filter error (MSE) to be lower with learning than with the ground truth parameters in the absence of learning, getting close to the performance of the optimal filter. 
This is shown in Fig.~\ref{Fig:5} in terms of trial-averaged learning curves. 
The normalized MSE with learning decreases within the time frame of $T=2000$ and converges below the MSE for the projection filter with fixed parameters set to the ground truth. 
The optimal performance was estimated by running a particle filter with prior importance function, resampling at every time-step, 1000 particles and parameters set to the ground truth \cite{Doucet:2000ui}. 

\section{Related approaches}
\label{sec:rel}
In this section, we attempt to review similar approaches for online maximum likelihood estimation, and their relations to our method.
We note that most of the literature on this topic is formulated for discrete-time systems, and we realize that the list of reviewed works is not exhaustive.
Some of the approaches for Hidden Markov Models (HMMs) discussed here are also surveyed in more detail in \cite{Cappe:2006wta} --\nocite{vanHandel:2008vz}\cite{Khreich:2012jn}.

\subsection{Recursive maximum-likelihood approaches}
\label{sec:discgrad}
This work is the continuous-time analogue of the online SGA algorithm of \cite{LeGland:1995gj,LeGland:1997kp} for HMMs.
The behavior of the algorithm is analyzed by casting it in the Robbins-Monro framework of stochastic approximations. 
We used a similar approach to studying convergence in Section~\ref{sec:conv}.
More recently, the convergence of discrete-time stochastic gradient algorithms for parameter estimation in HMMs was studied under more general conditions \cite{Tadic:2010ba}.
To our knowledge, it is an open problem to obtain a similarly general result for continuous-time models such as the one in this paper.
Regarding estimation in a discretized diffusion model, particle algorithms have been discussed in \cite{Poyiadjis:2006fb,Poyiadjis:2011dc}.

\subsection{Prediction error algorithms}
Another stochastic approximation scheme is the recursive minimum prediction error scheme (see \cite{LeGland:1997kp} and \cite{Collings:1994fc}) for HMMs.
Instead of finding maxima of the likelihood, it finds minima of the average (squared) prediction error, i.e. the error between the observations and the predicted observations.
In our continuous-time model, the prediction error is given by the infinitesimal pseudo-innovation increment $dY_t-\tilde h_t dt$. 
Formal differentiation of $(dY_t-\tilde h_t dt)^2$ with respect to the parameter yields the same parameter update rule as that derived in Section~\ref{sec:methods}.
While a rigorous analysis has not been done, it seems natural to conjecture that recursive maximum likelihood and recursive minimum prediction error are equivalent in continuous time.

\subsection{Online EM}
\label{sec:oEM}
Expectation maximization (EM) is a well-known method for offline parameter learning in partially observed stochastic systems \cite{Dempster:1977fi,Dembo:1986cv}.
It is based on the following application of Jensen's inequality:
\begin{equation}
\begin{split}
\mathcal{L}_t(\theta)-\mathcal{L}_t(\tilde\theta)&=\log\mathds{E}_{\tilde\theta}\brackets{\frac{dP_{\theta}}{dP_{\tilde\theta}}\cond{}\mathcal{F}_t^Y} \\
&\geq \mathds{E}_{\tilde\theta}\brackets{\log\frac{dP_{\theta}}{dP_{\tilde\theta}}\cond{}\mathcal{F}_t^Y} \doteq Q_t(\theta,\tilde\theta).
\end{split}
\end{equation}
Since $Q_t(\tilde\theta,\tilde\theta)=0$, by maximizing $Q_t(\theta,\tilde\theta)$ with respect to $\theta$ (for fixed $\tilde\theta$), we obtain a non-negative change in the likelihood.
EM thus produces a sequence of parameter estimates $\tilde\theta_k$, $k=0,1,2,...$ with non-decreasing likelihood by iterating the following procedure:
compute the quantity $Q_t(\theta,\tilde\theta_k)$ (the `expectation' or `E' step in EM), then set $\tilde\theta_{k+1}=\text{argmax}_{\theta} Q_t(\theta,\tilde\theta_k)$  (the `maximization' or `M' step in EM). 

If a parametrization is chosen such that the complete-data log-likelihood\footnote{We note that a limitation of EM in the continuous-time model is that the identification of parameters of the diffusion term $g_{\theta}$ has to be treated differently from that of drift parameters in $f_{\theta}$ and $h_{\theta}$. 
This is due to the fact that there is no reference measure for the complete model that is independent of the diffusion parameters. 
The parameters of the diffusion term are therefore not included in $\theta$, but are estimated separately from the quadratic variations of hidden state and observation.
This issue is discussed in more detail in \cite{Charalambous:2000jw}, Section IV-B. 
This issue is avoided in the gradient-based method here because the reference measure restricted to the observations is independent of all parameters, including the ones of the diffusion term.} takes the form of an exponential family, i.e. $\Psi(\theta)\cdot S_t$, where $\Psi$ is a vector-valued function of the parameters and $S_t$ is a vector of functionals of the hidden state and observation trajectories, then $Q_t(\theta,\tilde\theta)=\Psi(\theta)\cdot\hat S_t(\tilde\theta)+R(\tilde\theta)$, where
\begin{equation}
\hat S_t(\tilde\theta)=\mathds{E}_{\tilde\theta}\brackets{S_t\cond{}\mathcal{F}_t^Y},
\end{equation}
and $R(\tilde\theta)$ is independent of $\theta$.
The `M' step can be done explicitly if the equation $\partial_{\theta}\Psi(\theta)\cdot\hat S_t(\tilde\theta)=0$ has a unique closed-form solution.
Meanwhile, the `E' step consists of computing $\hat S_t(\tilde\theta)$, which involves certain nonlinear smoothed functionals of the forms 
\begin{align*}
&\mathds{E}_{\tilde\theta}\brackets{\int_0^t \varphi_1(X_s) dX_s\cond{}\mathcal{F}_t^Y},\\
&\mathds{E}_{\tilde\theta}\brackets{\int_0^t \varphi_2(X_s) dY_s\cond{}\mathcal{F}_t^Y},\\
&\mathds{E}_{\tilde\theta}\brackets{\int_0^t \varphi_3(X_s) ds\cond{}\mathcal{F}_t^Y},
\end{align*}
with possibly distinct integrands $\varphi_1,\varphi_2,\varphi_3$.
In general, these smoothed functionals are computed using a forward-backward smoothing algorithm, which is not suitable for online learning.
In a few select cases, the smoothed functionals admit a finite-dimensional solution (see \cite{Zeitouni:1986wi} and the remarks on p.99 of \cite{Dembo:1986cv}), or even a finite-dimensional \emph{recursive} solution (see \cite{Elliott:1997dl,Charalambous:2000jw,Charalambous:2003jn}).

In \cite{Charalambous:2000jw}, the smoothed functionals of the linear-Gaussian model are expressed (using the Fisher identity) in terms of derivatives of the incomplete-data log-likelihood, or a generalization thereof.
This enables a recursive computation of the smoothed quantities of interest, and the auxiliary variables that need to be integrated (called sensitivity equations) are very similar to Eqs.~\eqref{KBfilterderiv1}-\eqref{KBfilterderivlast}.
The relation between smoothed functionals and the sensitivity equations has been known for a long time (see \cite{Gupta:1974bw} and Section 10.2 in \cite{Cappe:2006wta}).

Several authors \cite{Krishnamurthy:1993gx}-\nocite{Mongillo:2008hq,Cappe:2009ik}\cite{Cappe:2011dw} have introduced the idea of a fully recursive form of EM, called \emph{online EM}.
In the references above, online EM has been explicitly formulated for HMMs and State Space Models (SSMs) by integrating the recursive smoothing algorithm using the online parameter estimate.
This stochastic approximation approach to EM is thus very similar to the gradient-based approach used here and in the references discussed in Section~\ref{sec:discgrad}.

In continuous-time diffusion models such as studied in this paper, the recursions found by \cite{Elliott:1997dl,Charalambous:2000jw,Charalambous:2003jn} can be directly applied if the model is linear.
We did this in order to do the comparison of SGA and online EM shown in Fig.~\ref{Fig:3}.
For this particular model, SGA and online EM are comparable in terms of computational complexity and rate of convergence.
In nonlinear models, online EM can be formulated by making use of recursive particle approximations of the smoothing functionals (e.g. by applying the methods in \cite{Cappe:2001ev,DelMoral:2010wv} to a suitable time discretization of the SDEs).
As an alternative, assumed-density or projection filters can be used to approximate the recursive smoothed functionals.
The full development of online EM in continuous time, as well as its convergence analysis, remains a topic for future research.

\subsection{State augmentation algorithms}
The idea is to treat the unknown parameter as a random variable that is either static ($d\theta_t = 0$) or has dynamics that are coupled to the hidden state.
In both cases the parameter may be estimated online by solving the filtering problem for the augmented state $(X_t,\theta_t)$.
While this presents clear advantages for known dynamics of the hidden parameter, it introduces a new parameter estimation problem for the parameters of the dynamics of $\theta_t$, called \emph{hyperparameters}. 
A static prior for $\theta_t$ is problematic because the resulting filter will usually not be stable, with negative implications (see \cite{Rebeschini:2015hs}) on the behavior of particle filters that are needed to solve the augmented filtering problem (but see \cite{Papavasiliou:2006km}, where stability conditions are discussed for the discrete-time case).
In addition, for many interesting models, the parameter space may be of much higher dimension than the state space, introducing high computational costs for filtering of the augmented state.

\subsection{Maximum-likelihood filtering and identification}
The opposite of state augmentation was explored in \cite{Moura:1986wi}, where the hidden state is also estimated via maximum likelihood, instead of the usual filtering paradigm using minimum mean-squared error.
Equations for the maximum-likelihood state and parameter estimates are then derived.
Although these equations are not directly suitable for recursive identification, they are very similar to the ones obtained by us in Section~\ref{sec:methods}.
It remains a curiosity that the approach of \cite{Moura:1986wi} has rarely been cited and has not been further developed.

\section{Conclusions}
\label{sec:disc}
The problem of estimating parameters in partially observed systems is old and relevant to many applications. 
However, the majority of the literature on this subject is written for discrete-time processes and for offline learning, while, despite of its enormous importance for filtering and control theory, the continuous-time case has received little attention.
Online gradient ascent in continuous time has only recently been studied in \cite{Sirignano:2017if}.
The use of a change of measure in order to express the likelihood function in terms of the filter is not new, but it seems to be underexploited. 
To the best of our knowledge, its only use in parameter estimation is in the technical report by \cite{Moura:1986wi}. 
We found it appropriate to revisit this approach and to extend the work of \cite{Sirignano:2017if} to the partially observable case.
Recently, the above results for the fully observed model have been strengthened to a central limit theorem for the parameter estimate, see \cite{Sirignano:2017wx}. 

The main difficulty and open problem is to find conditions on the generative model that are easy to verify, sufficient for the convergence of the algorithm, and not too restrictive.
Currently, the most promising avenue for obtaining such conditions is by settling open questions regarding the ergodicity of the approximate filter in terms of the exact one, and then using the general theory that guarantees ergodicity of the exact filter.
The latter is relatively easy to check compared to the explicit conditions on the approximate filter that we currently give in Section~\ref{sec:convB}.
We hope that these open questions will be addressed in the future.

Let us briefly comment on the numerical examples that we provided.
As we showed numerically, the algorithm is capable of improving filter performance even if the models are unidentifiable and the learning rate constant, even though this cannot be expected. 
In addition, the second numerical example showed that the performance of the filter can be improved even beyond what is possible with fixed parameters. 
This result could lead to new ways of improving the performance of approximate filters by using the additional degrees of freedom given by the online parameter estimates for both adaptation (learning) and reduced filter error. 
It remains to be explored whether this feature applies to a large enough class of approximate filters to be useful for practical applications.

As we showed also in comparison with the Online EM algorithm, these na\"ive methods exhibit rather slow convergence rates and cannot compete with fast offline methods such as second-order optimization methods or Nesterov's accelerated gradient.
However, the main aim of this article is to advance the theoretical understanding of convergence using continuous-time theory.
Based on this, it remains a topic for future research to study the convergence of more elaborate algorithms such as the ones mentioned above. 

\appendices
\section{Proof of proposition 1}
\label{propproof2}
\begin{enumerate}[label=(\roman*)]
\item We have
\begin{multline}
\frac{1}{t}\mathcal{L}_t(\theta)=\frac{1}{t}\int_0^t l(\mathcal{X}_s(\theta),\theta)ds\\
+\frac{1}{t}\int_0^t\psi_h(\theta,M_s(\theta))\cdot dV_s.
\label{prop1split}
\end{multline}
By Condition~\ref{condfinitedim}(i), the first term on the RHS converges to $\int_{\mathds{R}^D}l(x,\theta)\mu_{\theta}(dx)=\mathcal{\tilde L}(\theta)$ a.s. as $t\rightarrow\infty$.
Consider the local martingale $\mathfrak{M}_t=\int_0^t\psi_h(\theta,M_s(\theta))\cdot dV_s$.
From It\^o isometry, Condition~\ref{condregularity}, and Condition~\ref{condfinitedim}(v), it follows that for $t$ large enough,
\begin{equation}
\begin{split}
\mathds{E}&\brackets{\parenths{\int_0^t\psi_h(\theta,M_s(\theta))\cdot dV_s}^2} \\
&=\mathds{E}\brackets{\int_0^t\norm{\psi_h(\theta,M_s(\theta))}^2ds} \\
&\leq \mathds{E}\brackets{\int_0^tC(1+||M_s(\theta)||^q)ds} \\
&\leq \mathds{E}\brackets{\int_0^tC(1+||\mathcal{X}_s(\theta)||^q)ds} \\
&\leq Ct\parenths{1+\mathds{E}[\sup_{s\leq t}||\mathcal{X}_s(\theta)||^q]} \\
&\leq Ct(1+C'\sqrt{t}).
\end{split}
\end{equation}
In short, for $t$ large enough, we have $\text{Var}[\mathfrak{M}_t]\leq K t^{3/2}$ for some $K>0$.
Therefore, 
\begin{equation}
\text{Var}\brackets{\tfrac{1}{t}\mathfrak{M}_t}\leq K t^{-1/2}\rightarrow 0, \quad t\rightarrow\infty,
\end{equation}
which means that the second term on the RHS of Eq.~\eqref{prop1split} converges to zero in $L_2$.

Now consider the process $ \mathfrak{\tilde M}_t=\tfrac{1}{t}\mathfrak{M}_t+\int_0^t\tfrac{1}{s^2}\mathfrak{M}_sds$.
By It\^o's lemma, this process is the local martingale given by $\int_0^t\tfrac{1}{s}\psi_h(M_s(\theta),\theta)\cdot dV_s$.
By applying It\^o isometry, Condition~\ref{condfinitedim}(v) and \ref{condregularity}, we obtain
\begin{multline}
\sup_{t>0}\mathds{E}\brackets{||\mathfrak{\tilde M}_t||^2} \leq \int_0^{\infty}\frac{\mathds{E}\brackets{||\psi_h(M_s(\theta),\theta)||^q}}{s^2}ds \\
\leq K \int_0^{\infty}\frac{1}{s^2}\parenths{1+\mathds{E}\brackets{||\mathcal{X}_s(\theta)||^2}}ds<\infty. \\
\end{multline}
By the martingale convergence theorem, there is a finite random variable $\mathfrak{\tilde M}_{\infty}$ such that $\mathfrak{\tilde M}_t\to \mathfrak{\tilde M}_{\infty}$ a.s. and in $L_2$.
Therefore also $\tfrac{1}{t}\mathfrak{M}_t$ converges a.s.

\item We have that $\partial_{\theta}\mathcal{\tilde L}(\theta)=\lim_{t\rightarrow\infty}\frac{1}{t}\partial_{\theta}\mathcal{L}_t(\theta)$, if the derivative exists and the limit exists a.s. 
Due to Condition~\ref{condregularity}, the derivative
\begin{multline}
\frac{1}{t}\partial_{\theta}\mathcal{L}_t(\theta)=\frac{1}{t}\int_0^t \partial_{\theta}l(\mathcal{X}_s(\theta),\theta)ds\\
+\frac{1}{t}\int_0^t \partial_{\theta}\psi_h(M_s(\theta),\theta) dV_s \\
 = \frac{1}{t}\int_0^t F(\mathcal{X}_s(\theta),\theta)^{\top}ds+\frac{1}{t}\int_0^tdV^{\top}_sH(M_s(\theta),\theta) 
\end{multline}
exists.
This converges to $\int_{\mathds{R}^D}F(x,\theta)^{\top}\mu_{\theta}(dx)$ by an argument analogous to the one in (i).

The representation of $\mathcal{H}$ in terms of the invariant measure and its derivative follows from Conditions~\ref{condfinitedim}(iii) and \ref{condregularity}.

\item This follows from (i) and the fact that $F$ is in $\mathds{\bar G}$ (Condition~\ref{condregularity}).

%

\item By Condition~\ref{condregularity}, $q,K>0$ can be chosen such that the functions $l,F,\partial_{\theta}F,H$ grow at most as $K(1+||x||^q)$ for all $\theta\in\Theta$.
From this and the first part of the present Lemma, it follows that
\begin{equation}
\begin{split}
\mathcal{\tilde L}(\theta)&=\int_{\mathds{R}^D}l(x,\theta)\mu_{\theta}(dx) \\
&\leq K\int_{\mathds{R}^D}(1+||x||^q)\mu_{\theta}(dx)\leq K K_q.
\end{split}
\end{equation}
By a similar calculation, we have 
\begin{equation}
\norm{g(\theta)} \leq K K_q.
\end{equation}
For $\norm{\mathcal{H}(\theta)}$, observe that
\begin{equation}
\begin{split}
\norm{\mathcal{H}(\theta)}&\leq\norm{\int_{\mathds{R}^D}\partial_{\theta}F(x,\theta)\mu_{\theta}(dx)} \\
&\quad + \norm{\int_{\mathds{R}^D}F(x,\theta)\nu_{\theta}(dx)} \\
&\leq K K_q + \norm{\int_{\mathds{R}^D}F(x,\theta)\nu_{\theta}(dx)},
\end{split}
\end{equation}
where the first term on the RHS was treated in the same way as in the bound for $\mathcal{\tilde L}(\theta)$ and $\norm{g(\theta)}$.
For the second term, we observe that
\begin{multline}
\norm{\int_{\mathds{R}^D}F(x,\theta)\nu_{\theta}(dx)}^2 = \sum_{i,j=1}^p\parenths{\int_{\mathds{R}^D}F_i(x,\theta)\nu_{\theta,j}(dx)}^2 \\
\leq \sum_{i,j=1}^p\parenths{\int_{\mathds{R}^D}|F_i(x,\theta)|\,|\nu_{\theta,j}(dx)|}^2 \\
\leq \sum_{i,j=1}^p\parenths{\int_{\mathds{R}^D}||F(x,\theta)||\,|\nu_{\theta,j}(dx)|}^2 
\leq p^2K^2K_q'^2
\end{multline}
The claimed inequality~\eqref{eqn:lemma3b} then follows by setting $C=3KK_q+pK  K_q'$. \IEEEQEDclosed
\end{enumerate}

\section{Lemmas}
\label{lemmas}
Here, we adapt the Lemmas of \cite{Sirignano:2017if} to fit the present setting.
As in \cite{Sirignano:2017if}, the proofs of the lemmas require results from \cite{Pardoux:2003wk}, but in a slightly more general form than what was needed in \cite{Sirignano:2017if}.
Despite the strong similarities between our proofs and the proofs in \cite{Sirignano:2017if}, for the convenience of the reader we shall write them out in full detail and in the appropriate notation.

For the Lemmas 1-4 below, we assume that conditions~\ref{condfinitedim}-\ref{condlearningrate} hold and that the first exit time from $\Theta$ is infinite (see the proof of Theorem 1).
In addition, we define the following.
Let $\kappa,\lambda>0$ and define the $(P_{\theta_0},\mathcal{F}_t)$-stopping times $\sigma_0=0$ and $\sigma_k,\tau_k$, $k\in\mathds{N}$ as
\begin{align}
\tau_k&\doteq \inf\braces{t>\sigma_{k-1}:\Vert g(\tilde\theta_t)\Vert\geq\kappa}, \\
\begin{split}
\sigma_k&\doteq \sup \Big\{t>\tau_k:\tfrac{1}{2}\Vert g(\tilde\theta_{\tau_k})\Vert \leq \Vert g(\tilde\theta_{s})\Vert \\ 
&\quad \leq 2\Vert g(\tilde\theta_{\tau_k})\Vert,s\in[\tau_k,t]  \text{ and } \int_{\tau_k}^t \gamma_s ds\leq\lambda \Big\}
\end{split}
\end{align}

\begin{lemma}
\label{lemma1}
Let $\eta>0$ and define 
\begin{equation}
\Gamma_{k,\eta}\doteq\int_{\tau_k}^{\sigma_k+\eta}\gamma_s\parenths{F(\mathcal{\tilde X}_s,\tilde\theta_s)-g(\tilde\theta_s)^{\top}}ds
\end{equation}
Then, with probability one,
\begin{equation}
\lim_{k\rightarrow\infty}\norm{\Gamma_{k,\eta}}=0.
\end{equation}
\end{lemma}
\begin{IEEEproof}
Consider the function $G(x,\theta)=F(x,\theta)-g(\theta)^{\top}$. 
By definition, we have 
\begin{equation}
\int_{\mathds{R}^D}G(x,\theta)\mu_{\theta}(dx)=0,
\end{equation}
and by Condition~\ref{condregularity} we have that the components of $G(x,\cdot)$ are in $\mathds{\bar G}^D$.
Therefore, by Condition~\ref{condfinitedim}(iv), the Poisson equation
\begin{equation}
\mathcal{A_X}v(x,\theta)=G(x,\theta), \quad \int_{\mathds{R}^D}v(x,\theta)\mu_{\theta}(dx)=0
\end{equation}
has a unique twice differentiable solution with 
\begin{equation}
||v(x,\theta)||+||\partial_{\theta}v(x,\theta)||+||\partial_{\theta}^2v(x,\theta)||\leq K'(1+||x||^{q'}).
\label{eqn:vgrowth}
\end{equation}
Let $u(t,x,\theta)=\gamma_t v(x,\theta)$, and apply It\^o's lemma to each component of $u$:
\begin{multline}
u_i(\sigma,\mathcal{\tilde X}_{\sigma},\tilde\theta_{\sigma}) - u_i(\tau,\mathcal{\tilde X}_{\tau},\tilde\theta_{\tau})  = \int_{\tau}^{\sigma}\partial_s u_i(s,\mathcal{\tilde X}_{s},\tilde\theta_{s}) ds \\
+ \int_{\tau}^{\sigma}\mathcal{A_X} u_i(s,\mathcal{\tilde X}_{s},\tilde\theta_{s}) ds
+ \int_{\tau}^{\sigma}\mathcal{A}_{\theta} u_i(s,\mathcal{\tilde X}_{s},\tilde\theta_{s}) ds \\
+ \int_{\tau}^{\sigma}\gamma_s\text{tr}\brackets{\hat{\Sigma}(\mathcal{\tilde X}_{s},\tilde\theta_{s})H(\mathcal{\tilde X}_{s},\tilde\theta_{s})\partial_{\theta}^{\top}\partial_{x}u_i(s,\mathcal{\tilde X}_{s},\tilde\theta_{s})} ds \\
+ \int_{\tau}^{\sigma}\partial_{x}u_i(s,\mathcal{\tilde X}_{s},\tilde\theta_{s})\Sigma(\mathcal{\tilde X}_{s},\tilde\theta_{s})d\mathcal{B}_s \\
+ \int_{\tau}^{\sigma}\gamma_s\partial_{\theta}u_i(s,\mathcal{\tilde X}_{s},\tilde\theta_{s})H(\mathcal{\tilde X}_{s},\tilde\theta_{s})^{\top}dV_s,
\end{multline}
where $\mathcal{A_X}$ and $\mathcal{A}_{\theta}$ are the infinitesimal generators of the processes $\mathcal{X}_t$ and $\tilde\theta_t$, respectively, $\hat{\Sigma}(x,\theta)$ denotes the $(D\times n_y)$-matrix consisting of the rows $n'+1,n'+2,...,n'+n_y$ of the matrix $\Sigma(x,\theta)$, and $\partial_{\theta}^{\top}\partial_{x}u_k(s,x,\theta)_{ij}=\partial_{\theta_i}\partial_{x_j}u_k(s,x,\theta)$.
Using the Poisson equation and the previous identity, we obtain
\begin{equation}
\begin{split}
\Gamma_{k,\eta} &=\int_{\tau_k}^{\sigma_k+\eta}\gamma_s\parenths{F(\mathcal{\tilde X}_s,\tilde\theta_s)-g(\tilde\theta_s)}^{\top}ds \\
&= \int_{\tau_k}^{\sigma_k+\eta}\gamma_s G(\mathcal{\tilde X}_s,\tilde\theta_s)ds = \int_{\tau_k}^{\sigma_k+\eta}\gamma_s \mathcal{A_X}v(\mathcal{\tilde X}_s,\tilde\theta_s)ds \\
&= \int_{\tau_k}^{\sigma_k+\eta}\mathcal{A_X}u(s,\mathcal{\tilde X}_s,\tilde\theta_s)ds \\
&= \gamma_{\sigma_k+\eta} v(\mathcal{\tilde X}_{\sigma_k+\eta},\tilde\theta_{\sigma_k+\eta})-\gamma_{\tau_k} v(\mathcal{\tilde X}_{\tau_k},\tilde\theta_{\tau_k})\\
&\quad -\int_{\tau_k}^{\sigma_k+\eta}\dot\gamma_s v(\mathcal{\tilde X}_{s},\tilde\theta_{s}) ds -\int_{\tau_k}^{\sigma_k+\eta}\gamma_s \mathcal{A}_{\theta} v(\mathcal{\tilde X}_{s},\tilde\theta_{s}) ds \\
&\quad -\int_{\tau_k}^{\sigma_k+\eta}\gamma_s^2 \text{tr}\Big[\hat{\Sigma}(\mathcal{\tilde X}_{s},\tilde\theta_{s})H(\mathcal{\tilde X}_{s},\tilde\theta_{s}) \\
&\quad\times\partial_{\theta}^{\top}\partial_{x}\Big] v(\mathcal{\tilde X}_{s},\tilde\theta_{s})ds \\
&\quad -\int_{\tau_k}^{\sigma_k+\eta}\gamma_s \partial_{x}v(\mathcal{\tilde X}_{s},\tilde\theta_{s})\Sigma(\mathcal{\tilde X}_{s},\tilde\theta_{s})d\mathcal{B}_s \\
&\quad -\int_{\tau_k}^{\sigma_k+\eta}\gamma_s^2 \partial_{\theta}v(\mathcal{\tilde X}_{s},\tilde\theta_{s})H(\mathcal{\tilde X}_{s},\tilde\theta_{s})^{\top}dV_s.
\end{split}
\label{Gammaeq}
\end{equation}

Define
\begin{equation}
J_t^{(1)}\doteq \gamma_t\sup_{s\leq t}||v(\mathcal{\tilde X}_s,\tilde\theta_s)||.
\end{equation}
By using Condition~\ref{condfinitedim}, we have
\begin{equation}
\begin{split}
\mathds{E}\brackets{\parenths{J_t^{(1)}}^2} &= \mathds{E}\brackets{\gamma_t^2\sup_{s\leq t}||v(\mathcal{\tilde X}_s,\tilde\theta_s)||^2} \\
& \leq K\gamma_t^2\mathds{E}\brackets{1+\sup_{s\leq t}||\mathcal{\tilde X}_s||^q} \\
& = K\gamma_t^2\parenths{1+\mathds{E}\brackets{\sup_{s\leq t}||\mathcal{\tilde X}_s||^q}} \\
& \leq K K'\gamma_t^2(1+\sqrt{t}) \\
& \leq K''\gamma_t^2\sqrt{t},
 \end{split}
 \label{J1ineq}
\end{equation}
where the first two inequalities use Conditions~\ref{condfinitedim}(iv) and (v), respectively.
We choose an $r>0$ such that $\gamma_t^2 t^{1/2+2r}\rightarrow0$ for $t\rightarrow\infty$ (this is possible due to Condition~\ref{condlearningrate}), and we pick $T>0$ large enough such that $\gamma_t^2 t^{1/2+2r}\leq 1$ for $t\geq T$.
In addition, for each $0<\delta<r$ we define the event
$A_{t,\delta} \doteq \{J_t^{(1)}t^{r-\delta}\geq 1\}$.
For $t\geq T$,
\begin{equation}
\begin{split}
\mathds{P}(A_{t,\delta} ) &\leq \mathds{E}\brackets{J_t^{(1)}t^{r-\delta}} 
 \leq \mathds{E}\brackets{\parenths{J_t^{(1)}}^2}t^{2r-2\delta} \\
& \leq K''\gamma_t^2 t^{1/2+2r-2\delta} 
 \leq K'' t^{-2\delta},
 \end{split}
 \label{PAineq}
\end{equation}
where Eq.~\eqref{J1ineq} was used in the second inequality.\footnote{The first inequality in \eqref{PAineq} is elementary: For a nonnegative random variable $Y$ with law $p$, we have
$$\mathds{P}(Y\geq 1) = \int_1^{\infty} p(dy) \leq \int_1^{\infty} y p(dy) \leq \int_0^{\infty} y p(dy) = \mathds{E}(Y). $$}
We therefore have that
\begin{equation}
\sum_{n=1}^{\infty}\mathds{P}(A_{2^n,\delta} )  < \infty.
\end{equation}
By the Borel-Cantelli Lemma, only finitely many events $A_{2^n,\delta}$ can occur. 
Therefore, there is a random index $n_0$ such that $J_{2^n}^{(1)}2^{n(r-\delta)}\leq 1$ for all $n\geq n_0$. 
Alternatively, we can say that there is a finite positive random variable $\xi$ and a deterministic $n_1\in\mathds{N}$ such that 
\begin{equation}
J_{2^n}^{(1)}2^{n(r-\delta)}\leq \xi, \quad n\geq n_1
\end{equation}
(e.g. choose $\xi=\max \{\max_{1\leq n'\leq n_0}J_{2^{n'}}^{(1)}2^{n'(r-\delta)},1 \} $).
For $t\in[2^n,2^{n+1}]$ and $n\geq n_1$, we therefore have
\begin{equation}
\begin{split}
J_t^{(1)}&=\gamma_t \sup_{s\leq t}||v(\mathcal{\tilde X}_s,\tilde\theta_s)|| \leq \gamma_{2^n} \sup_{s\leq t}||v(\mathcal{\tilde X}_s,\tilde\theta_s)|| \\
&\leq \gamma_{2^n} \sup_{s\leq 2^{n+1}}||v(\mathcal{\tilde X}_s,\tilde\theta_s)|| \\
&\leq K\gamma_{2^{n+1}} \sup_{s\leq 2^{n+1}}||v(\mathcal{\tilde X}_s,\tilde\theta_s)|| \\
&= K J_{2^{n+1}}^{(1)} \leq K \frac{\xi}{2^{(n+1)(r-\delta)}} \leq K\frac{\xi}{t^{r-\delta}},
 \end{split}
\end{equation}
and as a consequence, $J_t^{(1)}\rightarrow 0$ a.s. as $t\rightarrow\infty$.

Next, define
\begin{equation}
\begin{split}
J_{t}^{(2)}&=\int_{0}^{t}\Big|\Big|\dot\gamma_s v(\mathcal{\tilde X}_{s},\tilde\theta_{s}) + \gamma_s \mathcal{A}_{\theta} v(\mathcal{\tilde X}_{s},\tilde\theta_{s}) \\
&\quad +\gamma_s^2 \text{tr}\brackets{\hat{\Sigma}(\mathcal{\tilde X}_{s},\tilde\theta_{s})H(\mathcal{\tilde X}_{s},\tilde\theta_{s})\partial_{\theta}^{\top}\partial_{x}} v(\mathcal{\tilde X}_{s},\tilde\theta_{s})\Big|\Big| ds \\
 \end{split}
\end{equation}
Due to the PGP of $H$, $\hat\Sigma$, and $v$ (Conditions~\ref{condfinitedim} and \ref{condregularity}), we have 
\begin{equation}
\begin{split}
\sup_{t>0}\mathds{E}\brackets{J_{t}^{(2)}}&\leq K \int_{0}^{\infty} \parenths{\dot\gamma_s+\gamma_s^2}\parenths{1+\mathds{E}[||\mathcal{\tilde X}_s||^q]}ds\\
&\leq KC \int_{0}^{\infty} \parenths{\dot\gamma_s+\gamma_s^2}ds < \infty.
 \end{split}
 \label{eqn:J2bound}
\end{equation}
In the first inequality we additionally used the fact that $\mathcal{A}_{\theta}$ contains at least a factor of $\gamma_t$, in the second one we relied on Condition~\ref{condfinitedim}(v) and in the third inequality we used Condition~\ref{condlearningrate}.
Thus $J_{t}^{(2)}$ converges to a finite random variable a.s.

Lastly, we have the term
\begin{equation}
\begin{split}
J_t^{(3)}&=\int_{0}^{t}\gamma_s \partial_{x}v(\mathcal{\tilde X}_{s},\tilde\theta_{s})\Sigma(\mathcal{\tilde X}_{s},\tilde\theta_{s})d\mathcal{B}_s \\
&\quad +\int_{0}^{t}\gamma_s^2 \partial_{\theta}v(\mathcal{\tilde X}_{s},\tilde\theta_{s})H(\mathcal{\tilde X}_{s},\tilde\theta_{s})^{\top}dV_s.
 \end{split}
\end{equation}
By using It\^o isometry and the same PGPs as in \eqref{eqn:J2bound}, we obtain
\begin{multline}
\sup_{t>0}\mathds{E}\brackets{||J_{t}^{(3)}||^2} = \int_0^{\infty}\gamma_s^2\mathds{E}\brackets{\norm{\partial_x v \Sigma}^2}ds \\
\quad + \int_0^{\infty}\gamma_s^4\mathds{E}\brackets{\norm{\partial_{\theta} v H^{\top}}^2}ds \\
\quad + 2\int_0^{\infty}\gamma_s^3\text{tr}\mathds{E}\brackets{\partial_{x}v\hat{\Sigma}H^{\top}\partial_{\theta}^{\top}v^{\top}}ds\\
\leq C K  \int_0^{\infty}\parenths{\gamma_s^2+\gamma_s^3+\gamma_s^4} \parenths{1+\mathds{E}[||\mathcal{\tilde X}_s||^q]} ds
\\
\leq C K C'  \int_0^{\infty}\parenths{\gamma_s^2+\gamma_s^3+\gamma_s^4}ds<\infty.
\end{multline}
Thus, by Doob's martingale convergence theorem, $J_{t}^{(3)}$ converges to a square integrable random variable a.s.

Finally, we note that
\begin{equation}
\begin{split}
||\Gamma_{k,\eta}|| &\leq J_{\sigma_k+\eta}^{(1)} + J_{\tau_k}^{(1)}+ J_{\sigma_k+\eta}^{(2)} - J_{\tau_k}^{(2)}\\
&\quad + ||J_{\sigma_k+\eta}^{(3)}-J_{\tau_k}^{(3)}|| \rightarrow 0, \quad k\rightarrow\infty,
 \end{split}
\end{equation}
which concludes the proof.
\end{IEEEproof}

\begin{lemma}
\label{lemma2}
Let $L$ be the Lipschitz constant of $g$.
Choose $\lambda>0$ such that for a given $\kappa>0$ (this is the parameter of the stopping times $\tau_k$) we have $3\lambda+\frac{\lambda}{4\kappa}=\frac{1}{2L}$.
For $k$ large enough and $\eta>0$ small enough, $\int_{\tau_k}^{\sigma_k+\eta}\gamma_s ds>\lambda$.
In addition, a.s., $\frac{\lambda}{2}\leq\int_{\tau_k}^{\sigma_k}\gamma_sds\leq\lambda$.
\end{lemma}
\begin{IEEEproof}
This proof goes through exactly like the proof of Lemma 3.2 in \cite{Sirignano:2017if}, with the only modification that the martingale in that proof takes the form
$$ \int_{0}^{t}\gamma_s \frac{g(\tilde\theta_s)}{R_s} H(\mathcal{\tilde X}_s,\tilde\theta_s)^{\top} dV_s. $$
\end{IEEEproof}

\begin{lemma}
\label{lemma4}
Suppose that $\tilde\theta_t\in\Theta$ for $t\geq 0$ and that there is an infinite number of intervals $[\tau_k,\sigma_k)$. 
There is a $\beta>0$ such that for $k>k_0$,
\begin{equation}
\mathcal{\tilde L}(\tilde\theta_{\sigma_k}) - \mathcal{\tilde L}(\tilde\theta_{\tau_k})\geq \beta
\end{equation}
a.s.
\end{lemma}
\begin{IEEEproof}
By using It\^o's lemma and the parameter update SDE \eqref{eqn:parameterSDE}, we obtain four terms:
\begin{multline}
\mathcal{\tilde L}(\tilde\theta_{\sigma_k}) - \mathcal{\tilde L}(\tilde\theta_{\tau_k}) = \int_{\tau_k}^{\sigma_k} \gamma_s\norm{g(\tilde\theta_s)}^2ds \\
+ \int_{\tau_k}^{\sigma_k}\gamma_s g(\tilde\theta_s) H(\mathcal{\tilde X}_s,\tilde\theta_s)^{\top} dV_s \\
+ \int_{\tau_k}^{\sigma_k}\frac{\gamma_s^2}{2}\text{tr}\brackets{H(\mathcal{\tilde X}_s,\tilde\theta_s)\mathcal{H}(\tilde\theta_s)H(\mathcal{\tilde X}_s,\tilde\theta_s)^{\top}}ds\\
+ \int_{\tau_k}^{\sigma_k}\gamma_s g(\tilde\theta_s) \cdot \brackets{F(\mathcal{\tilde X}_s,\tilde\theta_s)-g(\tilde\theta_s)^{\top}}ds \\
= \Omega_{1,k}+ \Omega_{2,k}+ \Omega_{3,k}+ \Omega_{4,k},
\end{multline}
where $\mathcal{H}$ is used to denote the Hessian of $\mathcal{\tilde L}$, see Eq.~\eqref{eqn:Hessian}.
By virtue of the definition of the stopping times and Lemmas~\ref{lemma1} and \ref{lemma2}, 
\begin{multline}
\Omega_{1,k}= \int_{\tau_k}^{\sigma_k} \gamma_s\norm{g(\tilde\theta_s)}^2ds \\ \geq \frac{\norm{g(\tilde\theta_{\tau_k})}^2}{4}\int_{\tau_k}^{\sigma_k} \gamma_s ds
\geq \frac{\norm{g(\tilde\theta_{\tau_k})}^2}{8}\lambda(\kappa).\\
\end{multline}
We define 
\begin{equation}
R_t=\begin{cases} 
||g(\tilde\theta_{\tau_k})||, & t\in[\tau_k,\sigma_k) \text{ for some k}\geq 1, \\
\kappa, & \text{else}
\end{cases},
\label{eqn:R}
\end{equation}
such that we can write
\begin{equation}
\begin{split}
\Omega_{2,k}&=\int_{\tau_k}^{\sigma_k}\gamma_s g(\tilde\theta_s) H(\mathcal{\tilde X}_s,\tilde\theta_s)^{\top} dV_s \\
&=\norm{g(\tilde\theta_{\tau_k})}\int_{\tau_k}^{\sigma_k}\gamma_s \frac{g(\tilde\theta_s)}{R_s} H(\mathcal{\tilde X}_s,\tilde\theta_s)^{\top} dV_s.
\end{split}
\label{eqn:Omega2}
\end{equation}
Since $||g(\tilde\theta_s)||/R_s\leq 2$, it follows from the It\^o isometry, Condition~\ref{condfinitedim}(v) and \ref{condregularity} that
\begin{multline}
\sup_{t\geq 0} \mathds{E}\brackets{\parenths{\int_0^{t} \gamma_s \frac{g(\tilde\theta_s)}{R_s} H(\mathcal{\tilde X}_s,\tilde\theta_s)^{\top} dV_s}^2}
\\\leq \sup_{t\geq 0}\int_0^{t} \mathds{E}\brackets{\gamma_s^2 \frac{||g(\tilde\theta_s)||^2}{R_s^2} ||H(\mathcal{\tilde X}_s,\tilde\theta_s)||^2}ds
\\\leq4\int_0^{\infty}\gamma_s^2\mathds{E}\brackets{\norm{H(\mathcal{\tilde X}_s,\tilde\theta_s)}^2}ds
\\\leq4K\int_0^{\infty}\gamma_s^2\parenths{1+\mathds{E}\brackets{||\mathcal{\tilde X}_s||^q}}ds<\infty.
\end{multline}
By Doob's martingale convergence theorem, the martingale $\mathfrak{M}_t=\int_0^{t} \gamma_s \frac{g(\tilde\theta_s)}{R_s} H(\mathcal{\tilde X}_s,\tilde\theta_s)^{\top} dV_s$ converges to a finite random variable $\mathfrak{M}$ as $t\rightarrow\infty$.
Thus for any $\epsilon>0$, there is a $k_0$ such that a.s. we have $\Omega_{2,k}\leq||g(\tilde\theta_{\tau_k})||\epsilon$ for all $k\geq k_0$.

Next, we consider $\Omega_{3,k}$. 
Using Conditions~\ref{condfinitedim} and \ref{condregularity} and Proposition~\ref{prop1}, we obtain
\begin{multline}
\sup_{t\geq 0} \mathds{E}\brackets{\abs{\int_0^{t}\frac{\gamma_s^2}{2}\text{tr}\brackets{H(\mathcal{\tilde X}_s,\tilde\theta_s)\mathcal{H}(\tilde\theta_s)H(\mathcal{\tilde X}_s,\tilde\theta_s)^{\top}}ds}}
\\ \leq \sup_{t\geq 0} \mathds{E}\brackets{\int_0^{t}\frac{\gamma_s^2}{2}\abs{\text{tr}\brackets{H(\mathcal{\tilde X}_s,\tilde\theta_s)\mathcal{H}(\tilde\theta_s)H(\mathcal{\tilde X}_s,\tilde\theta_s)^{\top}}}ds}
\\ \leq \int_0^{\infty} \frac{\gamma_s^2}{2} \mathds{E}\brackets{||H(\mathcal{\tilde X}_s,\tilde\theta_s)||^2||\mathcal{H}(\tilde\theta_s)||}ds
\\ \leq K \int_0^{\infty} \frac{\gamma_s^2}{2} \parenths{1+\mathds{E}\brackets{||\mathcal{\tilde X}_s||^{q}}}ds<\infty,
\end{multline}
from which it follows that $$\int_0^{t}\frac{\gamma_s^2}{2}\text{tr}\brackets{H(\mathcal{\tilde X}_s,\tilde\theta_s)\mathcal{H}(\tilde\theta_s)H(\mathcal{\tilde X}_s,\tilde\theta_s)^{\top}}ds$$ converges to a finite random variable as $t\rightarrow\infty$.
Thus, a.s., $\Omega_{3,k}$ must converge to zero as $k\rightarrow\infty$.

Finally, we consider the term $\Omega_{4,k}$ and define the function $G(x,\theta)=g(\theta) \cdot \brackets{F(x,\theta)-g(\theta)^{\top}}$, which by definition of $g$ satisfies $\int_{\mathds{R}^D}G(x,\theta)\mu_{\theta}(dx)=0$.
By Condition~\ref{condfinitedim}(iv), for each $\theta\in\Theta$ the Poisson equation $\mathcal{A_X}v(x,\theta)=G(x,\theta)$ (where $\mathcal{A_X}$ is the infinitesimal generator of the process $\mathcal{X}_t$) has a unique solution $v$ with $\int_{\mathds{R}^D}v(x,\theta)\mu_{\theta}(dx)=0$.
Let $u(t,x,\theta)\doteq\gamma_tv(x,\theta)$ and apply It\^o's lemma
\begin{multline}
u(\sigma,\mathcal{\tilde X}_{\sigma},\tilde\theta_{\sigma}) - u(\tau,\mathcal{\tilde X}_{\tau},\tilde\theta_{\tau})  = \int_{\tau}^{\sigma}\partial_s u(s,\mathcal{\tilde X}_{s},\tilde\theta_{s}) ds \\
+ \int_{\tau}^{\sigma}\mathcal{A_X} u(s,\mathcal{\tilde X}_{s},\tilde\theta_{s}) ds
+ \int_{\tau}^{\sigma}\mathcal{A}_{\theta} u(s,\mathcal{\tilde X}_{s},\tilde\theta_{s}) ds \\
+ \int_{\tau}^{\sigma}\gamma_s\text{tr}\brackets{\hat{\Sigma}(\mathcal{\tilde X}_{s},\tilde\theta_{s})H(\mathcal{\tilde X}_{s},\tilde\theta_{s})\partial_{x}\partial_{\theta}^{\top}u(s,\mathcal{\tilde X}_{s},\tilde\theta_{s})} ds \\
+ \int_{\tau}^{\sigma}\partial_{x}u(s,\mathcal{\tilde X}_{s},\tilde\theta_{s})\Sigma(\mathcal{\tilde X}_{s},\tilde\theta_{s})d\mathcal{B}_s \\
+ \int_{\tau}^{\sigma}\gamma_s\partial_{\theta}u(s,\mathcal{\tilde X}_{s},\tilde\theta_{s})H(\mathcal{\tilde X}_{s},\tilde\theta_{s})^{\top}dV_s,
\end{multline}
where $\hat{\Sigma}(x,\theta)$ denotes the $(D\times n_y)$-matrix consisting of the rows $n'+1,n'+2,...,n'+n_y$ of the matrix $\Sigma(x,\theta)$, and $\partial_{x}\partial_{\theta}^{\top}u(s,x,\theta)_{ij}=\partial_{\theta_i}\partial_{x_j}u(s,x,\theta)$.
Using the Poisson equation, we obtain
\begin{equation*}
\begin{split}
\Omega_{4,k}&=\int_{\tau_k}^{\sigma_k}\gamma_s g(\tilde\theta_s) \cdot \brackets{F(\mathcal{\tilde X}_s,\tilde\theta_s)-g(\tilde\theta_s)^{\top}}ds \\
&= \int_{\tau_k}^{\sigma_k}\gamma_s G(\mathcal{\tilde X}_s,\tilde\theta_s)ds = \int_{\tau_k}^{\sigma_k}\gamma_s \mathcal{A_X}v(\mathcal{\tilde X}_s,\tilde\theta_s)ds \\
&= \int_{\tau_k}^{\sigma_k}\mathcal{A_X}u(s,\mathcal{\tilde X}_s,\tilde\theta_s)ds,
\end{split}
\end{equation*}
which, by using the previous identity, turns into
\begin{equation}
\begin{split}
&= \gamma_{\sigma_k} v(\mathcal{\tilde X}_{\sigma_k},\tilde\theta_{\sigma_k})-\gamma_{\tau_k} v(\mathcal{\tilde X}_{\tau_k},\tilde\theta_{\tau_k})\\
&\quad -\int_{\tau_k}^{\sigma_k}\partial_s \gamma_s v(\mathcal{\tilde X}_{s},\tilde\theta_{s}) ds -\int_{\tau_k}^{\sigma_k}\gamma_s \mathcal{A}_{\theta} v(\mathcal{\tilde X}_{s},\tilde\theta_{s}) ds \\
&\quad -\int_{\tau_k}^{\sigma_k}\gamma_s^2 \text{tr}\brackets{\hat{\Sigma}(\mathcal{\tilde X}_{s},\tilde\theta_{s})H(\mathcal{\tilde X}_{s},\tilde\theta_{s})\partial_{\theta}^{\top}\partial_{x}v(\mathcal{\tilde X}_{s},\tilde\theta_{s})} ds \\
&\quad -\int_{\tau_k}^{\sigma_k}\gamma_s \partial_{x}v(\mathcal{\tilde X}_{s},\tilde\theta_{s})\Sigma(\mathcal{\tilde X}_{s},\tilde\theta_{s})d\mathcal{B}_s \\
&\quad -\int_{\tau_k}^{\sigma_k}\gamma_s^2 \partial_{\theta}v(\mathcal{\tilde X}_{s},\tilde\theta_{s})H(\mathcal{\tilde X}_{s},\tilde\theta_{s})^{\top}dV_s.
\end{split}
\end{equation}

By following the steps in the proof of Lemma 1, we find that $\Omega_{4,k}\rightarrow 0$ as $k\rightarrow\infty$ a.s.

For all $\epsilon>0$ a.s., we have for $k$ large enough
\begin{multline}
\mathcal{\tilde L}(\tilde\theta_{\sigma_k}) - \mathcal{\tilde L}(\tilde\theta_{\tau_k}) = \Omega_{1,k}+ \Omega_{2,k}+ \Omega_{3,k}+ \Omega_{4,k} \\
\geq \Omega_{1,k}- || \Omega_{2,k} || - || \Omega_{3,k} || - || \Omega_{4,k} || \\
\geq \frac{1}{8}\lambda(\kappa) || g(\tilde\theta_{\tau_k}) || ^2-\epsilon|| g(\tilde\theta_{\tau_k}) || -2\epsilon.
\end{multline}
The Lemma then follows by choosing $\epsilon=\min\{\frac{\lambda(\kappa)\kappa^2}{32},\frac{\lambda(\kappa)}{32}\}$ and $\beta=\frac{\lambda(\kappa)\kappa^2}{32}$.
\end{IEEEproof}

\begin{lemma}
\label{lemma5}
Under the conditions of Lemma \ref{lemma4} there is a $0<\beta_1<\beta$ such that for $k>k_0$
\begin{equation}
\mathcal{\tilde L}(\tilde\theta_{\tau_k}) - \mathcal{\tilde L}(\tilde\theta_{\sigma_{k-1}})\geq -\beta_1
\end{equation}
a.s.
\end{lemma}
\begin{IEEEproof}
As in Lemma~\ref{lemma4}, we obtain
\begin{multline}
\mathcal{\tilde L}(\tilde\theta_{\tau_k}) - \mathcal{\tilde L}(\tilde\theta_{\sigma_{k-1}}) \geq  \int_{\sigma_{k-1}}^{\tau_k}\gamma_s g(\tilde\theta_s) H(\mathcal{\tilde X}_s,\tilde\theta_s)^{\top} dV_s \\
+ \int_{\sigma_{k-1}}^{\tau_k}\frac{\gamma_s^2}{2}\text{tr}\brackets{H(\mathcal{\tilde X}_s,\tilde\theta_s)\mathcal{H}(\tilde\theta_s)H(\mathcal{\tilde X}_s,\tilde\theta_s)^{\top}}ds\\
+ \int_{\sigma_{k-1}}^{\tau_k}\gamma_s g(\tilde\theta_s) \cdot \brackets{F(\mathcal{\tilde X}_s,\tilde\theta_s)-g(\tilde\theta_s)^{\top}}ds.
\end{multline}
It is sufficient to show that the RHS converges to zero a.s.
Due to Eq.~\eqref{eqn:R}, the first term can be rewritten as
\begin{equation}
\kappa\int_{\sigma_{k-1}}^{\tau_k}\gamma_s \frac{g(\tilde\theta_s)}{R_s} H(\mathcal{\tilde X}_s,\tilde\theta_s)^{\top} dV_s.
\end{equation}
Using the argument from the proof of Lemma~\ref{lemma4}, this converges to zero a.s. as $k\rightarrow\infty$.
The treatment of the second and third terms is identical to the treatment of the terms $\Omega_{3,k}$ and $\Omega_{4,k}$ in the proof of Lemma~\ref{lemma4}.
\end{IEEEproof}

\section{Proof of Theorem 1}
\label{proof}
First, define the first exit time from $\Theta$
\begin{equation}
\tau=\inf\braces{t\geq 0:\; \tilde\theta_t\notin\Theta}.
\end{equation}
If $\tau<\infty$, since the paths of $\tilde\theta_t$ are continuous, we have $\tilde\theta_{\tau}\in\partial\Theta$. 
Furthermore, since $d\tilde\theta_t=0$ on $\partial\Theta$, we have $\tilde\theta_t\in\partial\Theta$ for all $t\geq\tau$.

Next, consider the case when $\tau=\infty$, which implies that $\tilde\theta_t\in\Theta$ for all $t\geq 0$.
Consider the case when there is a finite number of stopping times $\tau_k$.
Then, there is a finite $T$ such that $||g(\tilde\theta_t)||<\kappa$ for $t\geq T$. 
Therefore, since $\kappa$ can be chosen arbitrarily small, $\lim_{t\rightarrow\infty} ||g(\tilde\theta_t)||=0$.
Next, suppose that the number of stopping times $\tau_k$ is infinite. 
By Lemmas \ref{lemma4} and \ref{lemma5} there is a $k_0$ and constants $\beta>\beta_1>0$ such that for all $k\geq k_0$ a.s. 
\begin{align}
\mathcal{\tilde L}(\tilde\theta_{\sigma_k}) - \mathcal{\tilde L}(\tilde\theta_{\tau_k})&\geq \beta \\
\mathcal{\tilde L}(\tilde\theta_{\tau_k}) - \mathcal{\tilde L}(\tilde\theta_{\sigma_{k-1}})&\geq -\beta_1 > -\beta.
\end{align}
Thus, we have 
\begin{multline}
\mathcal{\tilde L}(\tilde\theta_{\tau_{n+1}})-\mathcal{\tilde L}(\tilde\theta_{\tau_{k_0}})\\
=\sum_{k=k_0}^{n}\brackets{\mathcal{\tilde L}(\tilde\theta_{\sigma_k}) - \mathcal{\tilde L}(\tilde\theta_{\tau_k})+\mathcal{\tilde L}(\tilde\theta_{\tau_{k+1}}) - \mathcal{\tilde L}(\tilde\theta_{\sigma_{k}})} \\
\geq (n+1-k_0)(\beta-\beta_1).
\end{multline}
Since $\beta-\beta_1>0$, when $n\rightarrow\infty$, $\mathcal{\tilde L}(\tilde\theta_{\tau_{n+1}})\rightarrow\infty$ a.s., and therefore $\mathcal{\tilde L}(\tilde\theta_t)\rightarrow\infty$ a.s.
This is in contradiction to Proposition~\ref{prop1}(iv), which states that $\mathcal{\tilde L}$ is bounded from above.
Therefore, there are a.s. only a finite number of stopping times $\tau_k$. \IEEEQEDclosed

\section*{Acknowledgments}
The authors would like to thank Philipp Harms, Ofer Zeitouni, Konstantinos Spiliopoulos, and Sean Meyn for helpful discussions.

\ifCLASSOPTIONcaptionsoff
  \newpage
\fi


\bibliographystyle{IEEEtran}
\bibliography{library}

%






\end{document}